\documentclass[a4paper, 10pt]{article}

\usepackage{amsmath}
\usepackage{amssymb}
\usepackage{latexsym}
\usepackage{theorem}
\usepackage[all]{xypic}

\newenvironment{proof}{\par \noindent{\bf Proof: }} {\nopagebreak \hspace{\stretch{1}} $\Box$ \par \mbox{}}
\newcommand{\noproof}{\hspace{\stretch{1}} $\Box$}
\newtheorem{theorem}{Theorem}[section]
\newtheorem{proposition}[theorem]{Proposition}
\newtheorem{lemma}[theorem]{Lemma}
\newtheorem{corollary}[theorem]{Corollary}
{\theorembodyfont{\rmfamily}
\newtheorem{definition}[theorem]{Definition}
\newtheorem{example}[theorem]{Example}
}
\newenvironment{theorem*}{\par \medskip \noindent{\bf Theorem }}{\par \mbox{}}
\newenvironment{lemma*}{\par \medskip \noindent{\bf Theorem }}{\par \mbox{}}

\newcommand{\Hom}{\mathop{Hom}}

\newcommand{\Ob}{\mathop{Ob}}

\newcommand{\gtp}{\mathop{\hat{\otimes}}}

\begin{document}

\title{$KK$-theory of $C^\ast$-categories and the analytic assembly map}

\author{Paul D. Mitchener \\
email: mitch@math.jussieu.fr} 

\maketitle

\section*{Abstract}
We define $KK$-theory spectra associated to $C^\ast$-categories and look at certain instances of the Kasparov product at this level.  This machinery is used to give a description of the analytic assembly map as a natural map of spectra.

Keywords: $C^\ast$-category, $KK$-theory, Assembly map

AMS 2000 Mathematics Subject Classification: 19K35, 55P42

\tableofcontents

\section{Introduction}

There are a number of different maps referred to as assembly maps, and in some cases several apparently different definitions exist of a given assembly map.  Such difficulties provide a motivation to explore how assembly maps can be characterised.  The following result is proved in \cite{WW}; an equivariant version of the same result exists in \cite{DL}.

\begin{theorem} \label{1}
Let $F$ be a homotopy-invariant functor from the category of topological spaces to the category of spectra.  Then there is a strongly excisive functor, $F^\%$, and a natural transformation, $\alpha \colon F^\%\rightarrow F$ such that the map $\alpha \colon F^\% (\textrm{point})\rightarrow F(\textrm{point})$ is a weak equivalence.

Further, the functor $F^\%$ and the natural transformation $\alpha$ are unique up to weak equivalence.  The map $\alpha \colon F^\%(X)\rightarrow F(X)$ is called the {\em assembly map} associated to the functor $F$.
\end{theorem}

Stripped of the technical language, this result says that a given assembly map which can be expressed as a map of spectra rather than simply as a group homomorphism is characterised by the fact that it maps something computable (this is what `excisive' means) to the object of interest in a natural way.

It is easy to apply this philosophy to a number of assembly maps, such as the assembly map in algebraic $K$-theory first appearing in \cite{Lod}, the assembly map in Waldhausen $K$-theory used to study the algebraic $K$-theory assembly map in \cite{BHM1} and \cite{BHM2}, the assembly map studied in \cite{Ti1}, and the topological assembly map used to formulate the Novikov conjecture in \cite{Ran}.  However, there are difficulties with using this approach to attempt to deal with the analytic assembly map 
$$\alpha_\star \colon K_\star (BG)\rightarrow K_\star (C^\star_\mathrm{max} G)$$
introduced to study the Novikov conjecture in \cite{Kas3}.

The problem is that we want to view the analytic assembly map as a special case of an assembly map
$$\alpha \colon F^\% (X)\rightarrow F(X)$$
in the sense of theorem \ref{1} defined for every space $X$.  Actually, it is not too hard to formulate the assembly map as a map of $K$-theory spectra
$$\alpha \colon {\mathbb K}_\mathrm{hom}(X)\rightarrow {\mathbb K} C^\star_\mathrm{max} \pi_1 (X)$$
where the spectrum on the left is related to the $K$-homology of the space $X$ and the spectrum on the right is related to the $K$-theory of the group $C^\ast$-algebra $C^\star_\mathrm{max} \pi_1 (X)$.

The main difficulty with this approach is that to define the fundamental group $\pi_1 (X)$ we need to choose a basepoint for the topological space $X$.  The assignment $X\mapsto {\mathbb K} C^\star_\mathrm{max} \pi_1 (X)$ is therefore not a functor on the category of topological spaces and all continuous maps, and theorem \ref{1} no longer applies.

We can get around this difficulty by replacing the fundamental group
$\pi_1 (X)$ by the fundamental groupoid $\pi (X)$, and looking at the
groupoid $C^\ast$-category $C^\star_\mathrm{max} \pi (X)$ introduced
in \cite{Mitch2}.  $K$-theory spectra associated to
$C^\ast$-categories are defined in \cite{Mitch2.5}, so we can form
the spectrum ${\mathbb K} C^\star_\mathrm{max} \pi (X)$, and the
assignment $X\mapsto {\mathbb K} C^\star_\mathrm{max} \pi (X)$ {\em
  is} a functor on the category of topological spaces and all
continuous maps.\footnote{The idea of using $C^\ast$-categories to
  avoid basepoint issues when defining spectrum-valued assembly maps
  first appeared in \cite{DL}.}  By theorem \ref{1} we obtain an assembly map
$$\alpha \colon {\mathbb K} ( C^\star_\mathrm{max} \pi )^\% (X)\rightarrow {\mathbb K} C^\star_\mathrm{max} \pi (X)$$

The main purpose of this article is to prove the following result.

\begin{theorem} \label{2}
The map
$$({\mathbb K} C^\star_\mathrm{max} \pi)^\% (X)\rightarrow \mathbb{K} C^\star_\mathrm{max} \pi (X)$$
is in fact the analytic assembly map.
\end{theorem}

By uniqueness of the assembly map associated to the functor $X\mapsto {\mathbb K} C^\star_\mathrm{max} \pi (X)$, all we need to do is show that the original assembly map can be written as a natural map of spectra
$$\alpha \colon {\mathbb K}_\mathrm{hom}(X)\rightarrow {\mathbb K} C^\star_\mathrm{max} \pi (X)$$

To write down the assembly map in this way, we need to have some notion of the $KK$-theory, and in particular $KK$-theory spectra, of $C^\ast$-categories.  The next three sections of this paper are devoted to developing the relevant features of such a theory.  We conclude in the final section by showing that the analytic assembly map can be written as a natural map of spectra in the way we desire.

\section*{Acknowledgements}
Much of this article is based on my D.~Phil dissertation.  I would therefore like to take this opportunity to thank my thesis advisors, John Roe and Ulrike Tillmann, for their frequent support.

\section{Preliminaries}

The purpose of this section is to summarise some of the concepts introduced elsewhere that are needed in the rest of this article.  In subsection \ref{cstarcat} we recall the definition of a $C^\ast$-category and in subsection \ref{lemon} we look at some concepts involved in defining $K$-theory.  In the final subsection, we prove a continuity result for the $K$-theory of $C^\ast$-categories that we will need later on.

\subsection{$C^\ast$-categories} \label{cstarcat}

Let $\mathbb F$ denote either the field of real numbers or the field of complex numbers.  Recall from \cite{Mi} that a category $\cal A$ is termed an {\em algebroid} if every morphism set is a vector space over the field $\mathbb F$, and composition of morphisms
$$\Hom (B,C)\times \Hom (A,B)\rightarrow \Hom (A,C)$$
is bilinear.  An {\em involution} on an algebroid $\cal A$ is a collection of maps $\Hom (A,B)\rightarrow \Hom (B,A)$, written $x\mapsto x^\star$, such that

\begin{itemize}

\item[$\bullet$] $(\alpha x + \beta y)^\star = \overline{\alpha} x^\star + \overline{\beta} y^\star$ for all scalars $\alpha , \beta \in {\mathbb F}$ and morphisms $x,y\in \Hom (A,B)$.

\item[$\bullet$] $(yx)^\star = x^\star y^\star$ for all morphisms $x\in \Hom (A,B)$ and $y\in \Hom (B,C)$

\item[$\bullet$] $(x^\star )^\star = x$ for all morphisms $x\in \Hom (A,B)$

\end{itemize}

We call an algebroid with involution a {\em $\star$-category}.  The definition of the objects on which we focus our attention in this article comes from \cite{GLR} and \cite{Mitch2}.

\begin{definition}
A {\em unital $C^\ast$-category} is a $\star$-category in which every morphism set $\Hom (A,B)$ is a Banach space, and

\begin{itemize}

\item[$\bullet$] The inequality $\| xy \| \leq  \| x \| \| y \|$ holds for all morphisms $x\in \Hom (B,C)$ and $y\in \Hom (A,B)$.

\item[$\bullet$] The {\em $C^\ast$-identity}, $\| x^\star x \| = \| x \|^2$, holds for all morphisms $x\in \Hom (A,B)$.

\item[$\bullet$] For every morphism $x\in \Hom (A,B)$, the product $x^\star x$ is a positive element of the $C^\ast$-algebra $\Hom (A,A)$.

\end{itemize}

A {\em non-unital $C^\ast$-category} is a collection of objects and morphisms similar to a $C^\ast$-category except that there need not exist identity morphisms $1\in \Hom (A,A)$.

\end{definition}

Note that a $C^\ast$-algebra can be considered to be a $C^\ast$-category with only one object.

We can form the category of all small unital
$C^\ast$-categories.\footnote{Recall that a category is said to be
  {\em small} if the class of all objects forms a set.  It is
  sometimes necessary to restrict our attention to small categories in
  order to avoid set-theoretic difficulties.}  The morphisms are {\em
  $C^\ast$-functors}, that is to say functors $F\colon {\cal
  A}\rightarrow {\cal B}$ between $C^\ast$-categories such that each
map $F\colon \Hom (A,B)\rightarrow \Hom (F(A),F(B))$ is linear, and
$F(x^\star ) = F(x)^\star$ for all morphisms $x$ in the category $\cal
A$.  We can similarly form the category of all small non-unital
$C^\ast$-categories.  Of course, in the definition of a
$C^\ast$-functor it is not necessary to assume that the
$C^\ast$-categories involved are small.

It is proved in \cite{Mitch2} that any $C^\ast$-functor is
automatically norm-decreasing, and therefore continuous.  Further, the
range of any $C^\ast$-functor is closed.

\begin{definition}
A $C^\ast$-category $\cal A$ is said to be {\em graded} if we can write each morphism set $\Hom (A,B)$ as a direct sum
$$\Hom (A,B) = \Hom (A,B)_0 \oplus \Hom (A,B)_1$$
of morphisms of degree $0$ and degree $1$ such that for composable morphisms $x$ and $y$ we have the formula
$$\deg (xy) = \deg (x) + \deg (y)$$

Here addition takes place modulo $2$.

A $C^\ast$-functor $F\colon {\cal A}\rightarrow {\cal B}$ between graded $C^\ast$-categories is termed a {\em graded $C^\ast$-functor} if
$$\deg (Fx) = \deg (x)$$
for every morphism $x$ in the category $\cal A$.
\end{definition}

We can consider an ungraded $C^\ast$-category to be equipped with the {\em trivial grading} defined by saying that every morphism is of degree $0$.  Our attitude is thus to view ungraded $C^\ast$-categories as special cases of graded $C^\ast$-categories.

There is a sensible notion of the spatial tensor product, ${\cal A}\gtp {\cal B}$, of graded $C^\ast$-categories $\cal A$ and $\cal B$.  The objects are pairs, written $A\otimes B$, for objects $A\in \Ob ({\cal A})$ and $B\in \Ob ({\cal B})$.  The morphism set $\Hom (A\otimes B,A'\otimes B')$ is a completion of the algebraic graded tensor product $\Hom (A,A')\ \hat{\odot}\ \Hom (B,B')$.  See section 7 of \cite{Mitch2} and definition 2.7 of \cite{Mitch2.5} for details.

\subsection{$K$-theory of $C^\ast$-categories} \label{lemon}

Before we are ready to look at the $K$-theory of $C^\ast$-categories, we need the following definition.

\begin{definition}
Let $\cal A$ be a small $C^\ast$-category.  Then we define the {\em additive completion}, ${\cal A}_\oplus$, to be the $\star$-category in which the objects are formal sums $A_1\oplus \cdots \oplus A_m$ where $A_i\in \Ob ({\cal A})$.  The morphism set $\Hom (A_1\oplus \cdots \oplus A_m , B_1\oplus \cdots \oplus B_n )$ is the set of matrices
$$\left\{ \left( \begin{array}{ccc}
x_{1,1} & \cdots & x_{1,m} \\
\vdots & \ddots & \vdots \\
x_{n,1} & \cdots & x_{n,m} \\
\end{array} \right)
\ |\ x_{ij} \in \Hom (A_j,B_i ) \right\}$$

Composition of morphisms is defined by matrix multiplication.  The involution is defined by the formula
$$\left( \begin{array}{ccc}
x_{1,1} & \cdots & x_{1,m} \\
\vdots & \ddots & \vdots \\
x_{n,1} & \cdots & x_{n,m} \\
\end{array} \right)^\star =
\left( \begin{array}{ccc}
x_{1,1}^\star & \cdots & x_{n,1}^\star \\
\vdots & \ddots & \vdots \\
x_{1,m}^\star & \cdots & x_{n,m}^\star \\
\end{array} \right)$$
\end{definition}

For any small $C^\ast$-category $\cal A$ there is an isometric $C^\ast$-functor into the category of all Hilbert spaces and bounded linear maps; see theorem 1.14 of \cite{GLR} and theorem 5.2 of \cite{Mitch2}.  Such a $C^\ast$-functor extends to a functor on the additive completion, ${\cal A}_\oplus$, and can be used to define a $C^\ast$-norm on the category ${\cal A}_\oplus$.  If the $C^\ast$-category $\cal A$ is graded then the additive completion ${\cal A}_\oplus$ is graded in the obvious way.

\begin{definition}
Let $\cal A$ be a small graded unital $C^\ast$-category.  Let $A\in \Ob ({\cal A}_\oplus)$.  Then we define the space of {\em supersymmetries}
$$SS(A) = \{ x\in \Hom (A,A)_1 \ |\ x=x^\star ,\ x^2 = 1 \}$$
\end{definition}

Given supersymmetries $x\in SS(A)$ and $y\in SS(B)$ we can form a supersymmetry
$$x\oplus y = \left( \begin{array}{cc}
x & 0 \\
0 & y \\
\end{array} \right) \in SS(A\oplus B)$$

Suppose that the spaces of supersymmetries $SS(A)$ are all non-empty.  Then we define $V_1({\cal A})$ to be the set 
$$\frac{\bigsqcup_{A\in \Ob ({\cal A}_\oplus )} \pi_0 SS(A)}{\langle x\oplus y\rangle \sim \langle y\oplus x\rangle }$$

Here we write $\langle x\rangle$ to denote the path-component of a supersymmetry $x\in SS(A)$.  The set $V_1({\cal A})$ is an Abelian semigroup with operation defined by the formula
$$\langle x \rangle + \langle y\rangle = \langle x\oplus y\rangle$$

\begin{definition}
Assume that the spaces of supersymmetries $SS(A)$ are all non-empty.  Then we define the group $G_1({\cal A})$ to be the Grothendieck completion of the semigroup $V_1({\cal A})$.  We define the $K$-theory group $K_1({\cal A})$ to be the set of formal differences
$$\{ \langle x\rangle -\langle y\rangle \in G_1({\cal A}) \ |\ x,y\in SS(A), A\in \Ob ({\cal A}_\oplus ) \}$$
\end{definition}

When $A$ is a graded $C^\ast$-algebra we recover from the above
definition the $K$-theory group $K_1 (A)$ defined in \cite{vD1}.

The main construction in \cite{Mitch2.5} is that of a functor,
$\mathbb K$, from the category of all small graded $C^\ast$-categories
to the category of symmetric spectra.\footnote{See \cite{HSS} for
  details of the theory of symmetric spectra.  The main feature of the
  theory for our purposes is that there is a well-behaved smash
  product in the category of symmetric spectra.}  When $\cal A$ is a
small unital graded $C^\ast$-category in which every space $SS(A)$ is
non-empty there is a natural isomorphism
$$K_1 ({\cal A}) \cong \pi_1 {\mathbb K}({\cal A})$$

For all other integers $n\in \mathbb Z$ we define the $K$-theory group $K_n ({\cal A})$ to be the stable homotopy group $\pi_n {\mathbb K}({\cal A})$.

The main elementary theorem in the $K$-theory of $C^\ast$-algebras is
the Bott periodicity theorem.  A version of the Bott periodicity
theorem is true for the $K$-theory of graded $C^\ast$-categories.  We
introduce it by first considering Clifford algebras, as was done for
graded $C^\ast$-algebras in \cite{vD2}.  See also \cite{ABS}, \cite{Kar1}, \cite{Kar2}, and \cite{Wo} for details on this approach to Bott periodicity.

\begin{definition}
Let $p$ and $q$ be natural numbers.  Then we define the {\em $(p,q)$-Clifford algebra}, $\mathbb F_{p,q}$, to be the algebra over the field $\mathbb F$ generated by elements
$$\{ e_1 , \ldots , e_p , f_1 , \ldots , f_q \}$$
that pairwise anti-commute and satisfy the formulae
$$e_i^2 = 1 \qquad f_j^2 = -1$$ 
\end{definition}

The Clifford algebra ${\mathbb F}_{p,q}$ is a graded $C^\ast$-algebra; the generators themselves are defined to be of degree $1$.  The following result, proved in \cite{Mitch2.5}, is the version of the Bott periodicity theorem that is most useful in the framework of this article.

\begin{theorem} \label{Bott}
Let $\cal A$ be a small graded $C^\ast$-category.  Then there is a natural stable equivalence of spectra
$$\Omega^q {\mathbb K}({\cal A}) \simeq \Omega^p {\mathbb K}({\cal A}\gtp \mathbb F_{p,q})$$
\noproof
\end{theorem}

We therefore have isomorphisms
$$K_{1-(p-q)}({\cal A}) = K_1 ({\cal A}\gtp {\mathbb F}_{p,q})$$

The above formulation of the Bott periodicity theorem enables us to
define the $K$-theory group $K_1 ({\cal A})$ without the assumption
that the spaces of supersymmetries $SS(A)$ are all non-empty.  To do
this, observe that the $C^\ast$-categories $\cal A$ and ${\cal
  A}\gtp {\mathbb F}_{1,1}$ have the same $K$-theory, and the
spaces of supersymmetries in the $C^\ast$-category ${\cal A}\gtp
{\mathbb F}_{1,1}$ are automatically non-empty.  In particular, the
$K$-theory of ungraded $C^\ast$-categories can be defined in this
way.  The usual $K$-theory groups of undgraded $C^\ast$-algebras are
recovered as a special case.

\begin{definition} \label{natiso}
Let $F, G\colon {\cal A}\rightarrow {\cal B}$ be graded $C^\ast$-functors
between unital graded $C^\ast$-categories.  Then a {\em natural isomorphism}
between $F$ and $G$ consists of a degree $0$ unitary morphism $U_A \in \Hom
(F(A),G(A))_{\cal B}$ for each object $A\in \Ob ({\cal A})$ such that
for every morphism $x\in \Hom (A,B)_{\cal A}$ the composites $U_B
F(x)$ and $F(x)U_A$ are equal.

A graded $C^\ast$-functor $F\colon {\cal A}\rightarrow {\cal B}$ between
unital $C^\ast$-categories is said to be an {\em equivalence} of
graded $C^\ast$-categories if there is a graded $C^\ast$-functor $G\colon {\cal B}\rightarrow {\cal A}$ such that the composites $FG$ and $GF$ are
naturally isomorphic to the identities $1_{\cal B}$ and $1_{\cal A}$ respectively.
\end{definition}

The following result is proved in \cite{Mitch2.5}

\begin{proposition} \label{natiso2}
Let $F\colon {\cal A}\rightarrow {\cal B}$ be an equivalence of small unital
graded $C^\ast$-categories.  Then the induced map $F_\star \colon {\mathbb
  K}({\cal A})\rightarrow {\mathbb K}({\cal B})$ is a stable
equivalence of $K$-theory spectra.
\noproof
\end{proposition}

In particular, if a small graded unital $C^\ast$-category is
equivalent to a $C^\ast$-algebra, it has the same $K$-theory.

\subsection{Continuity} \label{contsec}

In this subsection we investigate `continuity' properties of the $K$-theory groups $K_n ({\cal A})$, that is to say we investigate what happens when we look at direct limits.

We begin by looking at a `lifting' result for supersymmetries.

\begin{lemma} \label{ssdense}
Let $(A_i, \phi_{ij})$ be a directed family of graded unital $C^\ast$-algebras, with direct limit the graded $C^\ast$-algebra $A_\infty$, equipped with morphisms $\phi_i \colon A_i\rightarrow A_\infty$.

Let $x\in SS(A_\infty)$, $\varepsilon >0$.  Then there exists a
$C^\ast$-algebra $A_j$ and a point $y\in SS(A_j)$ such that
$$\| \phi_j (y) - x\| <  \varepsilon$$
\end{lemma}

\begin{proof}
It is a well-known fact from the theory of direct limits of $C^\ast$-algebras that the union $\cup_{i\in I} \phi_i (A_i)$ must be a dense subset of the $C^\ast$-algebra $A_\infty$; see for example \cite{W-O}, appendix L for details.  Choose a real number $\varepsilon_1 >0$.  Then there exists an element $a\in A_i$ such that $\| \phi_i (a) -x \| <  \varepsilon_1$,

Set $$b= \frac { a+ a^\star }{2}$$

Then the element $b$ is self-adjoint, and since $\phi_i$ is a morphism of $C^\ast$-algebras and $x$ is self-adjoint, we can see that $\| \phi_i (b) -x \| < \varepsilon_1$.  A similar trick gives us a degree $1$ self-adjoint element $c\in A_i$ such that
$$\| \phi_i (c) -x \| < \varepsilon_1$$

Now, without loss of generality, assume that $\varepsilon_1 < \frac{1}{3}$.  Observe
$$\begin{array}{rcl}
\| \phi_i (c)^2 -1 \| & \leq & \| \phi_i (c)^2 - \phi_i (c) x \| + \| \phi_i (c)x - x^2 \| \\
& \leq  & \| \phi_i(c) \| \| \phi_i(c) -x \| + \| x \| \| \phi_i (c) - x \| \\
& \leq & (\varepsilon_1 + 2\| x \| ) \| \phi_i (c) - x \| \\
& \leq & (\varepsilon_1 + 2) \varepsilon_1 \\
\end{array}$$

Define $\varepsilon_2 = \varepsilon_1 (2+\varepsilon_1)$.  Then $\varepsilon_1 < \frac{1}{3}$ so $\varepsilon_2<1$.  Therefore $\| \phi_i(c)^2 -1\| <1$.  Hence $\| \phi_i ( c^2 -1 )\| <1$.  But by definition of the norm on the direct limit $A_\infty$
$$\| \phi_i (c^2 -1) \| = \lim \sup \{ \| \phi_{ij} (c^2 -1) \| \ |\ j\geq i \}$$

So there exists an element $d = \phi_{ji}(c) \in A_j$ such that $ \| \phi_j (d) - x \| < \varepsilon_1$ and $\| d^2 -1 \| < \varepsilon_2$.

Since the element $d$ is self-adjoint, we have spectrum
$$\sigma (d) \subseteq (-1 -\varepsilon_3 , -1 + \varepsilon_3) \cup (1-\varepsilon_3, 1+\varepsilon_3)$$
where $\varepsilon_3 = \varepsilon_2 ^{1/2}$.  Certainly, $\varepsilon_3 <1$ so we may define a continuous function $f \colon \sigma (d)\rightarrow \mathbb R$ by writing
$$f(t) = \left\{ {1\ \ t\geq 0 \atop -1\ \ t\leq 0} \right.$$

Hence, by functional calculus there is an self-adjoint element $y = f(d)\in A_j$ such that $y^2 =1$.  By the formula expressing the norm of a self-adjoint element of a $C^\ast$-algebra in terms of its spectrum (see for example chapter 12 of \cite{Rud} for details), we have the inequality
$$\begin{array}{rcl}
\| y - d \| & = & \| f(d) -d \| \\
& \leq & \sup_{t\in (-1 -\varepsilon_3 , -1 + \varepsilon_3) \cup (1-\varepsilon_3, 1+\varepsilon_3) } | f(t) - t | \\
& < & \varepsilon_3 \\
\end{array}$$

To summarise, $y\in SS(A_j)$ and
$$\begin{array}{rcl}
\| \phi_j (y) - x \| & \leq & \| \phi_j \| \| y-c \| + \| \phi_j (c) -x \| \\
& < & \varepsilon_3 + \varepsilon_1 \\
\end{array}$$

Recall that $\varepsilon_3 = (\varepsilon _1 (2+\varepsilon_1))^{1/2}$.  Given $\varepsilon >0$, we can easily find a real number $\varepsilon_1 > 0$ such that $\varepsilon_3 + \varepsilon_1 < \varepsilon$.  Hence there is a point $y\in SS(A_j)$ such that $\| \phi_j(y) - x \| <\varepsilon $.
\end{proof}

It is possible to show that direct limits always exist in the category
of small graded $C^\ast$-categories.  The proof very closely follows
the proof of the existence of direct limits in the category of
$C^\ast$-algebras.  The interested reader can consult \cite{Mitch} for
details.  

\begin{theorem}
Let $({\cal A}_i, \phi_{ij})$ be a directed family of small unital graded $C^\ast$-categories with direct limit ${\cal A}_\infty$.  Then the family of groups $(K_1({\cal A}_i), {\phi_{ij}}_\star )$ has direct limit $K_1({\cal A}_\infty)$.
\end{theorem}

\begin{proof}
Observe that the additive completion $({\cal A}_\infty)_\oplus$ is the
direct limit of the family of additive completions $(({\cal
  A}_i)_\oplus, \phi_{ij})$.  Let the $C^\ast$-category ${\cal
  A}_\infty$ be equipped with $C^\ast$-functors $\phi_i \colon {\cal
  A}_i\rightarrow {\cal A}_\infty$ such that $\phi_j \phi_{ij} =
\phi_i$ whenever $i\leq j$.  Without loss of generality, assume that
for all objects $A_i \in \Ob ({\cal A}_i )$ the set of supersymmetries
$SS (A_i)$ is not empty.\footnote{As we have already mentioned, we
  can, if necessary, replace each $C^\ast$-category ${\cal A}_i$ by
  the tensor product ${\cal A}_i\gtp{\mathbb F}_{1,1}$ without
  affecting the $K$-theory.}

Consider an object $A\in \Ob (({\cal A}_\infty)_\oplus )$.  By the universal property of direct limits we can find an object $A_i \in \Ob (({\cal A}_i)_\oplus )$ such that $\phi_i (A_i) = A$.  Whenever $i\leq j$ define an object $A_j = \phi_{ij} (A_i)$.  Then we have a directed family of $C^\ast$-algebras $(\Hom (A_j ,A_j ) , \phi_{ij})_{j\geq i}$ with direct limit $\Hom (A,A)$.

Recall from proposition 1.11 of \cite{Mitch2.5} that if supersymmetries $x,y\in SS(A)$ are sufficiently close then the path-components $\langle x\rangle$ and $\langle y\rangle$ are equal.  Hence, by lemma \ref{ssdense}, for any supersymmetry $x\in SS(A)$ there is a supersymmetry $x_i \in \Hom (A_i ,A_i )$ such that ${\phi_i}_\star \langle x_i \rangle = \langle x\rangle$.  It follows that we can write
$$K_1 ({\cal A}) = \bigcup_{i\in I} {\phi_i}_\star K_1 ({\cal A}_i)$$

Suppose that $G$ is a group equipped with maps $\psi_i \colon K_1 ({\cal A}_i)\rightarrow G$ such that $\psi_j {\phi_{ij}}_\star = \psi_i$ whenever $i\leq j$.  Then we can define a homomorphism $\beta \colon K_1 ({\cal A}_\infty)\rightarrow G$ by writing $\beta {\phi_i}_\star (\alpha ) = \psi_i (\alpha )$ whenever $\alpha \in K_1 ({\cal A}_i)$.  Thus the group $K_1 ({\cal A}_\infty)$ has the universal property and so must be the direct limit of the family $(K_1({\cal A}_i) , {\phi_{ij}}_\star )$.
\end{proof}

If $\cal A$ is a non-unital $C^\ast$-category then there are unital $C^\ast$-categories ${\cal A}^+$ and ${\mathbb F}_{\cal A}$ together with a natural quotient $C^\ast$-functor $\pi \colon {\cal A}^+\rightarrow {\mathbb F}_{\cal A}$.  The $K$-theory group $K_1 ({\cal A})$ can be defined to be the kernel
$$K_1 ({\cal A}) = \ker (\pi_\star \colon K_1 ({\cal A}^+) \rightarrow K_1 ({\mathbb F}_{\cal A}))$$

See section 2 of \cite{Mitch2} and proposition 2.26 of \cite{Mitch2.5} for further details.

\begin{corollary} \label{nudl}
Let $({\cal A}_i, \phi_{ij})$ be a directed family of non-unital
graded $C^\ast$-categories with direct limit ${\cal A}_\infty$.  Then
the directed family of groups $(K_1({\cal A}_i), {\phi_{ij}}_\star )$ has direct limit $K_1({\cal A}_\infty)$.
\end{corollary}

\begin{proof}
We have short exact sequences of groups
$$0\rightarrow K_1({\cal A}_i)\rightarrow K_1({\cal A}_i^+)\rightarrow K_1(\mathbb F_{{\cal A}_i})\rightarrow 0$$
and the maps $\phi_{ij}$ induce commutative diagrams
$$\begin{array}{ccccc}
K_1({\cal A}_i) & \rightarrow & K_1({\cal A}_i^+) & \rightarrow & K_1(\mathbb F_{{\cal A}_i}) \\
\downarrow & & \downarrow & & \downarrow \\
K_1({\cal A}_j) & \rightarrow & K_1({\cal A}_j^+) & \rightarrow & K_1(\mathbb F_{{\cal A}_j}) \\
\end{array}$$

By the above theorem we know that the families of groups $(K_1({\cal A}_i^+))$ and $(K_1(\mathbb F_{{\cal A}_i}))$ have direct limits $K_1({\cal A}_\infty ^+)$ and $K_1(\mathbb F_{{\cal A}_\infty})$ respectively.  Let $G$ be the direct limit of the family of groups $(K_1({\cal A}_i))$.

It is a well-known result of homological algebra (see for example \cite{Sp} or \cite{We}) that the direct limit of a family of short exact sequences is a short exact sequence.  Thus we have a short exact sequence
$$0\rightarrow G\rightarrow K_1({\cal A}_\infty^+)\rightarrow K_1(\mathbb F_{{\cal A}_\infty})\rightarrow 0$$

By definition of the group $K_1({\cal A}_\infty)$ it follows that $G\cong K_1({\cal A}_\infty )$.
\end{proof}

\begin{corollary} \label{kcont}
Let $({\cal A}_i, \phi_{ij})$ be a directed family of graded
$C^\ast$-categories with direct limit ${\cal A}_\infty$.  Then the
directed family of groups $(K_n({\cal A}_i), {\phi_{ij}}_\star )$ has direct limit $K_n({\cal A}_\infty)$.
\end{corollary}

\begin{proof}
Observe that the family $({\cal A}_i\gtp {\mathbb F}_{p,q}, \phi_{ij}\otimes 1)$ has direct limit ${\cal A}_\infty \gtp {\mathbb F}_{p,q}$.  The desired conclusion now follows from the above result and the Bott periodicity theorem.
\end{proof}

As an application of the above continuity theorem we prove a stability result that will be important to us later on.

\begin{definition}
We define the algebra of {\em compact operators}, $\cal K$, to be the $C^\ast$-algebra direct limit of matrix algebras
$${\cal K} = \lim_{\longrightarrow \atop k} M_k({\mathbb F})$$
under the inclusions $M_k({\mathbb F})\rightarrow M_{k+1}({\mathbb F})$ defined by the formula
$$x\mapsto \left( \begin{array}{cc}
x & 0 \\
0 & 0 \\ \end{array} \right)$$

We equip the $C^\ast$-algebra $\cal K$ with the trivial grading.
\end{definition}

The following result is easy to see.

\begin{lemma} \label{tklem}
Let $\cal A$ be any $C^\ast$-category.  Then
$${\cal A}\gtp {\cal K} = \lim_{\longrightarrow \atop k} ({\cal A}\gtp M_k(\mathbb F))$$
\noproof
\end{lemma}

Let us write $p$ for the element of the $C^\ast$-algebra $\cal K$ arising from the element $1\in M_1({\mathbb F})$.\footnote{Or more generally for any rank $1$ projection in the $C^\ast$-algebra $\cal K$}  
Then our previous machinery gives us the following result.

\begin{theorem} \label{bigstab}
The $C^\ast$-functor ${\cal A}\rightarrow {\cal A}\gtp {\cal K}$ defined by writing $a\mapsto a\gtp p$ induces an isomorphism
$$K_n({\cal A})\cong K_n({\cal A}\gtp {\cal K})$$
at the level of $K$-theory groups.
\end{theorem}

\begin{proof}
It is easy to check that each induced map $K_n(M_k({\cal A}))\rightarrow K_n(M_{k+1}({\cal A}))$ is an isomorphism.  The result now follows by corollary \ref{kcont} and lemma \ref{tklem}.
\end{proof}

The $C^\ast$-functor ${\cal A}\rightarrow {\cal A}\gtp {\cal K}$ therefore induces a stable equivalence of symmetric spectra ${\mathbb K}({\cal A})\rightarrow {\mathbb K}({\cal A}\gtp {\cal K})$

\section{Hilbert modules over $C^\ast$-categories}

The $KK$-theory of $C^\ast$-algebras is defined in terms of Hilbert modules.  The same is true for the $KK$-theory of $C^\ast$-categories.  We look at the relevant concepts in this section.

\subsection{Hilbert Modules} \label{modules}

Let $\cal A$ be a $C^\ast$-category.  Then a {\em right $\cal A$-module} is a linear contravariant functor, $\cal E$, from the category $\cal A$ to the category of vector spaces.  We use the notation
$$\eta x = {\cal E}(x){\cal E}(\eta )$$
to denote the action of a morphism $x\in \Hom (A,B)_{\cal A}$ on a vector $\eta \in {\cal E}(A)$.  It is similarly possible to define {\em left ${\cal A}$-modules}.

A right $\cal A$-module, $\cal E$, is said to be {\em countably
  generated} if there is a countable set
$$\Omega \subseteq \bigcup_{A\in \Ob ({\cal A})} {\cal E}(A)$$
such that for each object $A\in \Ob ({\cal A})$, every element of the vector
space ${\cal E}(A)$ is a finite linear combination of elements of the
form $\eta x$, where $x\in \Hom (A,B)$ and $\eta \in \Omega \cap {\cal E}(B)$.

Recall from \cite{Mitch2} that a {\em semi-inner product} on a right $\cal A$-module $\cal E$ is a collection of maps $\langle -,-\rangle \colon {\cal E}(B)\times {\cal E}(A)\rightarrow \Hom (A,B)_{\cal A}$ such that

\begin{itemize}

\item[$\bullet$] For all vectors $\eta \in {\cal E}(B)$, $\xi , \zeta \in {\cal E}(C)$, and morphisms $x,y \in \Hom(A,C)$ we have the formula
$$\langle \eta , \xi x + \zeta y \rangle = \langle \eta , \xi \rangle x + \langle \eta , \zeta \rangle y$$

\item[$\bullet$] $\langle \eta , \xi \rangle = \langle \xi ,\eta \rangle^\star$

\item[$\bullet$] For any vector $\eta \in {\cal E}(A)$ the product $\langle \eta ,\eta \rangle$ is a positive element of the $C^\ast$-algebra $\Hom (A,A)$.

\end{itemize}

A semi-inner product is called an {\em inner product} if for any vector $\eta \in {\cal E}(A)$ the product $\langle \eta ,\eta \rangle$ is zero if and only if $\eta =0$.  A {\em Hilbert $\cal A$-module} is a right $\cal A$-module, $\cal E$, equipped with an inner product such that each space ${\cal E}(A)$ is complete with respect to the norm
$$\| x \| = \| \langle x ,x\rangle \|^\frac{1}{2}$$

\begin{definition}
We call a Hilbert $\cal A$-module $\cal E$ {\em
  countably generated} if there is a countably generated right $\cal A$-module ${\cal E}_0$ such that the space ${\cal E}_0 (A)$ is a dense
  subset of the space ${\cal E}(A)$ for every object $A\in \Ob ({\cal
  B})$.

The countable set $\Omega$ which generates the right $\cal A$-module
${\cal E}_0$ is referred to as a {\em generating set} for the Hilbert
$\cal A$-module $\cal E$.
\end{definition}

The above definition of a countably generated differs from the
definition given in \cite{Mitch2} but agrees with that of
\cite{Jo2}.  If the $C^\ast$-category $\cal A$ is unital and equivalent to a
$C^\ast$-algebra, the above definition is the same as that in
\cite{Mitch2}.

\begin{example} \label{oneeq1}
Let $\cal A$ be a $C^\ast$-category.
Then for each object $A\in \Ob ({\cal A})$ we have a Hilbert $\cal A$-module $\Hom (-,A)_{\cal A}$ with spaces $\Hom (C,A)_{\cal A}$.
The $\cal A$-action is defined by composition of morphisms, and the
inner product is defined by the formula
$$\langle x,y \rangle = x^\star y$$
\end{example}

\begin{example} \label{standardmodule}
Let $\cal A$ be a $C^\ast$-category.  For each pair of objects $C,A\in
\Ob ({\cal A})$, let us write ${\cal H}(C,A)_{\cal A}$ to denote the
space of sequences, $(x_n )$, such that each element $x_n$ belongs to
the morphism set $\Hom (C,A)_{\cal A}$, and the series
$$\sum_n x_n^\star x_n$$
converges in norm in the space $\Hom (C,C)_{\cal A}$.

Then we have a Hilbert $\cal A$-module ${\cal H}(-,A)_{\cal A}$
defined by associating the space ${\cal H}(C,A)_{\cal A}$ to the
object $C$.  The action of the $C^\ast$-category $\cal A$ is defined
by the formula
$$( x_n ) y = (x_n y)$$

The inner product is defined by the formula
$$\langle (x_n ), (y_n ) \rangle = \sum_n x_n^\star y_n$$

We can show that the inner product is well-defined by using the
Cauchy-Schwarz inequality.
\end{example}

We need to have a notion of the direct sum of a collection of Hilbert
$\cal A$-modules.

\begin{definition}
Let $\{ {\cal E}_{\lambda} \ |\ \lambda \in \Lambda \}$ be a countable
collection of Hilbert $\cal A$-modules.  Write
$$\prod_{\lambda \in \Lambda}{\cal E}_\lambda (A)$$
to denote the space of sequences of vectors $(\eta_\lambda )$ such
that $\eta_\lambda \in {\cal E}_\lambda (A)$ and the series
$$\sum_{\lambda \in \Lambda} \langle \eta_\lambda , \eta_\lambda
\rangle$$
converges in the $C^\ast$-algebra $\Hom (A,A)_{\cal A}$.

Then we define the {\em direct sum}
$$\prod_{\lambda \in \Lambda}{\cal E}_\lambda$$ 
to be the Hilbert $\cal A$-module associating the vector space
$\prod_{\lambda \in \Lambda}{\cal E}_\lambda (A)$ to the object $A\in
\Ob ({\cal A})$.
\end{definition}

The action of the category $\cal A$ on the direct sum $\prod_{\lambda
  \in \Lambda}{\cal E}_\lambda$ is defined in
the obvious way.  The inner product is defined by the formula
$$\langle (\eta_\lambda ), (\xi_\lambda ) \rangle = \sum_{\lambda \in
  \Lambda} \langle \eta_\lambda , \xi_\lambda \rangle$$

In particular, one can form the direct sum, ${\cal E}\oplus {\cal E}'$
of two Hilbert $\cal A$-modules $\cal E$ and ${\cal E}'$.  If $\{
{\cal E}_{\lambda} \ |\ \lambda \in \Lambda \}$ is a countable
collection of countably generated Hilbert $\cal A$-modules, the direct sum
$\prod_{\lambda \in \Lambda}{\cal E}_\lambda (A)$ is countably generated. 

Recall that a $C^\ast$-algebra $A$ is said to be {\em $\sigma$-unital}
if there is a countable approximate unit, that is a countable set
$$\{ e_\lambda \ |\ \lambda \in \Lambda \}$$
such that
$$\inf_{\lambda \in \Lambda} \| xe_\lambda - x \| =0 \qquad
\inf_{\lambda \in \Lambda} \| e_\lambda x - x \| =0$$

\begin{definition}
A $C^\ast$-category $\cal A$ is called {\em $\sigma$-unital} if each
endomorphism set $\Hom (A,A)_{\cal A}$ is a $\sigma$-unital
$C^\ast$-algebra.
\end{definition}

If $\cal A$ is a $\sigma$-unital $C^\ast$-category, the Hilbert $\cal
B$-modules $\Hom (-,A)_{\cal A}$ and ${\cal H}(-,A)_{\cal A}$ are
countably generated.

Suppose that $\cal A$ and $\cal B$ are $C^\ast$-categories, that $\cal E$ is a Hilbert $\cal A$-module, and that $\cal F$ is a Hilbert $\cal B$-module.  Then we define the {\em outer tensor product}, ${\cal E}\otimes {\cal F}$, to be the Hilbert ${\cal A}\otimes {\cal B}$-module in which the space $({\cal E}\otimes {\cal F})(A\otimes B)$ is the completion of the algebraic tensor product ${\cal E}(A)\odot {\cal F}(B)$ with respect to the norm defined by the inner product
$$\langle \eta_1 \otimes \xi_1 , \eta_2 \otimes \xi_2 \rangle = \langle \eta_1 ,\eta_2 \rangle \otimes \langle \xi_1 ,\xi_2 \rangle$$

It is proved in \cite{Mitch2} that the outer tensor product ${\cal E}\otimes {\cal F}$ is a Hilbert ${\cal A}\otimes {\cal B}$-module.

\begin{definition}
Let $\cal A$ be a graded $C^\ast$-category.  Then a {\em Hilbert $\cal A$-module} $\cal E$ is called {\em graded} if it admits decompositions ${\cal E}(A) = {\cal E}(A)_0 \oplus {\cal E}(A)_1$ into vectors of degree $0$ and vectors of degree $1$ such that

\begin{itemize}

\item[$\bullet$] $\deg (\eta x) = \deg (\eta) + \deg (x)$ for all vectors $\eta \in {\cal E}(B)$ and morphisms $x\in \Hom (A,B)_{\cal A}$

\item[$\bullet$] $\deg (\langle \eta , \xi \rangle ) = \deg (\eta ) + \deg (\xi )$ for all vectors $\eta \in {\cal E}(B)$ and $\xi \in {\cal E}(A)$.

\end{itemize}

Here all addition takes place modulo $2$.
\end{definition}

The direct sum and outer tensor product of graded Hilbert modules can
be graded in the obvious way.  If $\cal A$ is a graded
$C^\ast$-category, the Hilbert $\cal A$-modules $\Hom (-,A)_{\cal A}$
and ${\cal H}(-,A)_{\cal A}$ inherit a grading from the
$C^\ast$-category $\cal A$.

Let $\cal A$ be a graded $C^\ast$-category, and let $\cal E$ and $\cal F$ be Hilbert $\cal A$-modules.  Then an {\em operator} $T\colon {\cal E}\rightarrow {\cal F}$ is a collection of maps $T_A \colon {\cal E}(A)\rightarrow {\cal F}(A)$ such that there exist maps $T_A^\star \colon {\cal F}(A)\rightarrow {\cal E}(A)$ satisfying the formula
$$\langle \eta , T_A \xi \rangle = \langle T_B^\star \eta , \xi \rangle$$
for all vectors $\eta \in {\cal F}(B)$ and $\xi \in {\cal E}(A)$.

If $T\colon {\cal E}\rightarrow {\cal F}$ is an operator it can be shown (see \cite{Mitch2}) that

\begin{itemize}

\item[$\bullet$] $T_A(\eta x + \xi y) = (T_B\eta )x + (T_B\xi )y$ for all vectors $\eta , \xi \in {\cal E}(B)$ and all morphisms $x,y\in \Hom (A,B)$.  In particular, the operator $T$ is natural.

\item[$\bullet$] Each map $T_A\colon {\cal E}(A)\rightarrow {\cal F}(A)$ is linear and continuous.

\item[$\bullet$] The collection of maps $T_A^\star$ is uniquely determined by the operator $T$, and defines an operator $T^\star$ such that $( T^\star )^\star = T$.

\end{itemize}

Although each map $T_A\colon {\cal E}(A)\rightarrow {\cal E}'(A)$ is continuous, the norm
$$\| T \| = \sup \{ \| T_A \| \ |\ A\in \Ob ({\cal A}) \}$$
need not be finite.  When the norm $\| T\|$ is finite we call the operator $T$ {\em bounded}.  We write ${\cal L}({\cal E},{\cal F})$ to denote the space of all bounded operators $T\colon {\cal E}\rightarrow {\cal F}$.

Let ${\cal L}({\cal A})$ be the category of all countably generated graded Hilbert $\cal A$-modules and bounded operators between them.  It is proved in \cite{Mitch2} that the category ${\cal L}({\cal A})$ is a $C^\ast$-category.  It can be graded by defining the degree of a bounded operator $T\colon {\cal E}\rightarrow {\cal E}'$ to be $0$ if $\deg (T\eta ) = \deg (\eta )$ for all vectors $\eta \in {\cal E}(A)$ and to be $1$ if $\deg (T\eta ) = \deg (\eta )+1$ for all vectors $\eta \in {\cal E}(A)$.

\begin{definition}
Let $\cal A$ be a $C^\ast$-category, and let $\cal E$ and $\cal F$ be
Hilbert $\cal A$-modules.  Then a {\em rank one operator} $T\colon
{\cal E}\rightarrow {\cal F}$ is an operator of the form
$$\zeta \mapsto \eta \langle \xi , \zeta \rangle$$
for elements $\eta \in {\cal F}$ and $\xi \in {\cal E}$.  We write
this operator $\eta \langle \xi , - \rangle$.  A {\em
  compact operator} is a norm-limit of finite linear combinations of rank one operators.
\end{definition}

The above definition of a compact operator differs from the
definition given in \cite{Mitch2} but agrees with that of
\cite{Jo2}.  If the $C^\ast$-category $\cal A$ is unital and equivalent to a
$C^\ast$-algebra, the above definition is the same as that in
\cite{Mitch2}.

We write ${\cal K}({\cal E},{\cal E}')$ to denote the space of compact operators from a graded Hilbert ${\cal A}$-module $\cal E$ to a graded Hilbert $\cal A$-module ${\cal E}'$.  We write ${\cal K}({\cal A})$ to denote the collection of all countably generated graded Hilbert $\cal A$-modules and compact operators. 
Recall from \cite{GLR} that a {\em $C^\ast$-ideal}, $\cal J$, in a $C^\ast$-category $\cal A$ is a $C^\ast$-subcategory such that the composite of a morphism in the category $\cal J$ and a morphism in the category $\cal A$ belongs to the category $\cal J$.  We can form the {\em quotient}, ${\cal A}/{\cal J}$, of a $C^\ast$-category by a $C^\ast$-ideal.

\begin{proposition} \label{kideall}
The collection ${\cal K}({\cal A})$ is a $C^\ast$-ideal in the $C^\ast$-category ${\cal L}({\cal A})$.
\end{proposition}

\begin{proof}
It is clear that each space ${\cal K}({\cal E},{\cal E}')$ is a closed subspace of the Banach space ${\cal L}({\cal E},{\cal E}')$.  It is easy to see that the composition of a finite rank operator and a bounded operator is of finite rank.  Hence the composition of a compact operator with any other bounded operator is compact.
\end{proof}

We write ${\cal Q}({\cal A})$ to stand for the quotient ${\cal L}({\cal A})/{\cal K}({\cal A})$

The final result we need in this section is a version of the Kasparov
stabilisation theorem for Hilbert modules over $C^\ast$-categories.
It is more general than the version proved in \cite{Mitch2}; such a
generalisation is needed to perform some of the constructions in
$KK$-theory to be found in this article.

\begin{theorem} \label{kst}
Let $\cal A$ be a $\sigma$-unital $C^\ast$-category, and let $\cal E$ be a countably
generated Hilbert $\cal A$-module.  Then there is a countable
collection, $\Gamma$, of objects of the $C^\ast$-category such that we
have an isomorphism of Hilbert $\cal A$-modules
$${\cal E}\oplus \left( \prod_{A\in \Gamma} {\cal H}(-,A)_{\cal A}
\right) \cong \prod_{A\in \Gamma} {\cal H}(-,A )_{\cal A}$$
\end{theorem}

Before launching ourselves into a proof of this result, we need a
lemma.  To formulate the lemma, let $\cal E$ and ${\cal E}'$ be Hilbert $\cal
B$-modules.  Let us say an operator $T\in {\cal L}({\cal E},{\cal
  E}')$ has {\em dense range} if for each object $A\in \Ob ({\cal A})$
the set
$$\{ T\eta \ |\ \eta \in {\cal E}(A) \}$$
is a dense subset of the space ${\cal E}'(A)$.

\begin{lemma} \label{Kaslem}
Let $\cal E$ and ${\cal E}'$ be Hilbert $\cal A$-modules.  Suppose we
have an operator $T\in {\cal L}({\cal E},{\cal E}')$ such that the
operator $T$ and its adjoint $T^\star$ both have dense range.  Then
the Hilbert $\cal A$-modules $\cal E$ and ${\cal E}'$ are isomorphic.
\end{lemma}

\begin{proof}
The desired conclusion follows immediately from the corresponding
result for Hilbert modules over $C^\ast$-algebras.  For a proof of the
result for Hilbert modules over $C^\ast$-algebras, see for example
proposition 3.8 of \cite{Lan}.
\end{proof}

Our proof now parallels the proof of the Kasparov stabilisation theorem
for Hilbert modules over $C^\ast$-algebras given in \cite{MP}. 

\mbox{}

\noindent
{\bf Proof of theorem \ref{kst}: }
Let 
$$\Omega \subseteq \bigcup_{A\in \Ob ({\cal A})} {\cal E}(A)$$
be a countable generating set for the Hilbert $\cal A$-module $\cal
E$.  Write 
$$\Gamma = \{ A\in \Ob ({\cal A}) \ |\ {\cal E}(A)\cap \Omega \neq
\emptyset \}$$
and let $(\eta_n )$ be a sequence of vectors in the set $\Omega$ such
that every element of $\Omega$ occurs infinitely often.

Since the $C^\ast$-category $\cal A$ is $\sigma$-unital we can
construct a sequence, $(e_n )$, of unit vectors in the union
$\cup_{A\in \Gamma} {\cal H}(A,A)_{\cal A}$
such that:

\begin{itemize}

\item If $\eta_n \in {\cal E}(A)$, then $e_n \in {\cal H}(A,A)_{\cal
    B}$

\item The set 
$$\{ e_n \ |\ \eta_n \in {\cal E}(A) \}$$
is an orthonormal basis for the space ${\cal H}(A,A)$.

\end{itemize}

Observe that the set of vectors $e_n$ may be considered a generating
set for the Hilbert $\cal A$-module
$$\prod_{A\in \Gamma} {\cal H}(-,A)_{\cal A}$$

The Hilbert $\cal A$-modules $\cal E$ and ${\cal H}(-,A)_{\cal A}$ can
be regarded as submodules of the direct sum
$${\cal E}\oplus \left( \prod_{A\in \Gamma} {\cal H}(-,A)_{\cal A} \right)$$
in the obvious way.  We can define a compact operator $T\in {\cal
  K}( \prod_{A\in \Gamma} {\cal H}(-,A)_{\cal A} , {\cal E}\oplus (
\prod_{A\in \Gamma} {\cal H}(-,A)_{\cal A} ))$ by the formula
$$T\eta = \sum_{n=1}^\infty \left( 2^{-n} \eta_n \langle e_n ,\eta
  \rangle + 4^{-n} e_n \langle e_n , \eta \rangle \right)$$

In particular:
$$Te_n = 2^{-n} \eta_n + 4^{-n} e_n$$

For every natural number $m$ such that $\eta_m = \eta_n$ we have the
formula
$$T(2^m e_m ) = \eta_n + 2^{-m} e_m$$

But by construction of the sequence $(\eta_n )$ there are infinitely
many such natural numbers $m$.  Hence the vector $\eta_n$ belongs to
the closure of the image of the operator $T$.  Further,
$$T(4^n e_n ) = 2^n \eta_n + e_n$$
so the vector $e_n$ must also be in the closure of the range of the
operator $T$.  But the set of vectors $e_n$ and $\eta_n$ is a
generating set for the Hilbert $\cal A$-module
$${\cal E}\oplus \left( \prod_{A\in \Gamma} {\cal H}(-,A)_{\cal A} \right)$$
so the operator $T$ has dense range.

Now,
$$T^\star \xi = \sum_{n=1}^\infty \left( 2^{-n} e_n \langle \eta_n ,\xi
  \rangle + 4^{-n} e_n \langle e_n , \xi \rangle \right)$$

Hence $T^\star (4^n e_n ) = e_n$ and the operator $T^\star$ also has
dense range.  Therefore, by lemma \ref{Kaslem}, the Hilbert $\cal
B$-modules
$${\cal E}\oplus \left( \prod_{A\in \Gamma} {\cal H}(-,A)_{\cal A} \right)$$
and
$$\prod_{A\in \Gamma} {\cal H}(-,A)_{\cal A}$$
are isomorphic and we are done.
\noproof

There is an obvious version of the Kasparov stabilisation theorem in
the graded case.

\subsection{Bimodules}

The following definition is of great relevance to the $KK$-theory of $C^\ast$-categories.

\begin{definition}
Let $\cal A$ and $\cal B$ be graded $C^\ast$-categories.  A graded {\em Hilbert $({\cal A},{\cal B})$-bimodule} is a graded $C^\ast$-functor ${\cal F} \colon {\cal A}\rightarrow {\cal L}({\cal B})$
\end{definition}

For each object $A\in \Ob ({\cal A})$ we write ${\cal F}(-,A)$ to denote the corresponding Hilbert $\cal B$-module.  Given another object $B\in \Ob ({\cal B})$ we write ${\cal F}(A,B)$ to denote the associated vector space.  For each morphism $x\in \Hom (A,A')_{\cal A}$ we write $x\colon {\cal F}(-,A)\rightarrow {\cal F}(-,A')$ to denote the induced operator ${\cal F}(x)$.

Thus a graded Hilbert $({\cal A},{\cal B})$-bimodule $\cal F$ consists of a collection of Hilbert $\cal B$-modules
$$\{ {\cal F}(-,A)\ |\ A\in \Ob ({\cal A}) \}$$
together with a bounded operator $x\colon {\cal F}(-,A)\rightarrow {\cal F}(-,A')$ for each morphism $x\in \Hom ({\cal A},{\cal B})_{\cal A}$.

\begin{example} \label{onedim}
Suppose that $\cal B$ is a $\sigma$-unital graded $C^\ast$-category.  Let $F\colon {\cal A}\rightarrow {\cal B}$ be a graded $C^\ast$-functor.
Recall from example \ref{oneeq1} that for each object $B\in \Ob ({\cal
  B})$ we have a graded Hilbert $\cal B$-module $\Hom (-,B)_{\cal B}$
with spaces $\Hom (C,B)_{\cal B}$.  The Hilbert $\cal B$-module $\Hom
(-,B)_{\cal B}$ is countably generated since the $C^\ast$-category $\cal B$ is $\sigma$-unital.

The category $\cal B$ can thus itself be considered a graded Hilbert $({\cal
  A},{\cal B})$-bimodule; the $C^\ast$-functor $F\colon {\cal
  A}\rightarrow {\cal B}$ associates the graded Hilbert $\cal B$-module $\Hom
  (-,F(A))_{\cal B}$ to the object $A\in \Ob ({\cal A})$.  Given a
  morphism $x\in \Hom (A,A')_{\cal A}$ we have an operator $F(x)
  \colon \Hom (-,F(A))_{\cal B}\rightarrow \Hom (-,F(A'))_{\cal B}$
  defined by composition in the category $\cal B$.
\end{example}

It is possible to form the tensor product of a Hilbert $\cal A$-module and a Hilbert $({\cal A}, {\cal B})$-bimodule.  Before we are ready to construct such a product we need some technical results.  Our approach here closely follows the analysis of tensor products of Hilbert modules over $C^\ast$-algebras made in \cite{Lan}.

\begin{lemma} \label{pretens1}
Let $\cal E$ be a Hilbert $\cal A$-module, and consider a bounded operator $T\in {\cal L}({\cal E},{\cal E})$.  Then the operator $T$ is a positive element of the $C^\ast$-algebra ${\cal L}({\cal E},{\cal E})$ if and only if the product $\langle \eta, T \eta \rangle$ is positive for all vectors $\eta \in {\cal E}(A)$.
\end{lemma}

\begin{proof}
Suppose that the operator $T$ is a positive element of the $C^\ast$-algebra ${\cal L}({\cal E},{\cal E})$.  Then we can write $T = S^\star S$ for some operator $S\in {\cal L}({\cal E},{\cal E})$.  Hence for any vector $\eta \in {\cal E}(A)$
$$\langle \eta, T\eta \rangle = \langle \eta , S^\star S\eta \rangle = \langle S\eta , S\eta \rangle$$
so the product $\langle \eta , T\eta \rangle$ is positive.

Conversely, suppose that the product $\langle \eta , T\eta \rangle$ is positive for all vectors $\eta \in {\cal E}(A)$.  Then by polarisation the operator $T$ is self-adjoint.  We can therefore write $T = R-S$ for positive operators $R,S\in {\cal L}({\cal E},{\cal E})$ such that $RS = 0$.  Let $\eta \in {\cal E}(A)$.  By hypothesis the product
$$\langle S\eta , T (S\eta )\rangle = \langle S\eta , -S^2 \eta \rangle = \langle \eta ,-S^3 \eta \rangle$$
is positive.  But the operator $S$ is positive so we must have $S = 0$.  Hence the operator $T$ is positive.
\end{proof}

\begin{lemma} \label{pretens2}
Let $\cal E$ be a Hilbert $\cal A$-module.  Suppose that we are given vectors $\eta_1 \in {\cal E}(A_1), \ldots , \eta_n \in {\cal E}(A_n)$.  Let $A$ be the $C^\ast$-algebra $\Hom (A_1\oplus \cdots \oplus A_n,A_1\oplus \cdots \oplus A_n)$.  Form the matrix
$$X = \left( \begin{array}{ccc}
\langle \eta_1 ,\eta_1 \rangle & \cdots & \langle \eta_1 , \eta_n \rangle \\
\vdots & \ddots & \vdots \\
\langle \eta_n ,\eta_1 \rangle & \cdots & \langle \eta_n , \eta_n \rangle \\
\end{array} \right) \in A$$

Then the matrix $X$ is a positive element of the $C^\ast$-algebra $A$.
\end{lemma}

\begin{proof}
Consider the Hilbert $A$-module ${\cal F} = \Hom (A_1,A_1)\oplus \cdots \oplus \Hom (A_n,A_n)$.  There is an injective morphism of $C^\ast$-algebras $\rho \colon A\rightarrow {\cal L}({\cal F},{\cal F})$ defined by allowing a matrix $T\in A$ to act on a vector $b\in {\cal F}$ by matrix multiplication.

Form the vector $\eta = (\eta_1 ,\ldots , \eta_n )\in {\cal E}(A_1)\oplus \cdots \oplus {\cal E}(A_n)$.  Then for any vector $b = (b_1 ,\ldots ,b_n )\in {\cal F}$ we know that
$$\langle b , Xb \rangle = \sum_{i,j} b_i^\star \langle \eta_i ,\eta_j \rangle b_j = \langle \eta b,\eta b \rangle$$
so the product $\langle b,Xb \rangle$ is positive.  Hence the matrix $X$ is positive by the previous lemma.
\end{proof}

A proof of our final preliminary lemma can be found in \cite{Lan}.

\begin{lemma} \label{pretens3}
Let $\cal E$ be a Hilbert module over a $C^\ast$-algebra $A$, and suppose that we are given a vector $\eta \in {\cal E}$ and a real number $\alpha \in (0,1)$.  Then there is an element $\eta'\in {\cal E}$ such that $\eta = \eta'\langle \eta ,\eta \rangle^\frac{\alpha}{2}$
\noproof
\end{lemma}

Let $\cal E$ be a Hilbert $\cal A$-module and let $\cal F$ be a Hilbert $({\cal A},{\cal B})$-bimodule.  Then for each object $B\in \Ob ({\cal B})$ the functor ${\cal F}(B,-)$ defined by sending the object $A\in \Ob ({\cal A})$ to the vector space ${\cal F}(B,A)$ is a left $\cal A$-module. 

\begin{definition}
We define the {\em algebraic tensor product} ${\cal E}\odot_{\cal A} {\cal F}(B,-)$ to be the vector space
$$\{ \lambda_1 (\eta_1 ,\xi_1 ) + \cdots + \lambda_n (\eta_n ,\xi_n) \ |\ \lambda_i \in {\mathbb F}, \eta_i \in {\cal E}(A) , \xi_i \in {\cal F}(A,B) , A\in \Ob ({\cal A}) \} / \sim$$
where $\sim$ is the equivalence relation defined by writing

\begin{itemize}

\item[$\bullet$] $(\eta_1 + \eta_2 , \xi ) \sim (\eta_1 , \xi ) + (\eta_2 ,\xi )$

\item[$\bullet$] $(\eta , \xi_1 + \xi_2 ) \sim (\eta , \xi_1 ) + (\eta , \xi_2 )$

\item[$\bullet$] $(\eta , F(x) \xi ) \sim (\eta x , \xi )$

\end{itemize}

\end{definition}

We write $\eta \otimes \xi$ for the equivalence class of the pair $(\eta , \xi )$.  A more category-theoretic description of the algebraic tensor product of a right module and a left module over a $C^\ast$-category can be found in \cite{Mitch2}.

The above algebraic tensor product can be graded in the obvious way.  When the grading is taken into account, we write it ${\cal E}\ \hat{\odot}_{\cal A}\ {\cal F}(B,-)$.

\begin{lemma}
The collection of vector spaces
$${\cal E}\odot_{\cal A}{\cal F} = \{ {\cal E}\odot_{\cal A}{\cal F}(B,-) \ |\ B\in \Ob (B) \}$$
is a right $\cal B$-module.  The action of the category $\cal B$ is defined by the formula
$$(\eta\otimes \xi)(x) = \eta \otimes (\xi x)$$

There is a semi-inner product defined by the formula
$$\langle \eta\otimes \xi ,\eta'\otimes \xi' \rangle  = \langle \xi , \langle \eta ,\eta' \rangle \xi' \rangle$$
\end{lemma}

\begin{proof}
The only non-trivial part to check is the fact that for any vector
$$\zeta = \eta_1 \otimes \xi_1 + \cdots + \eta_n \otimes \xi_n \in {\cal E}\odot_{\cal A}{\cal F}(B)$$ 
where $\eta_i \in {\cal E}(A_i)$ and $\xi_j \in {\cal F}(B,A_i)$ the product $\langle \zeta ,\zeta \rangle$ is positive.  Observe that
$$\langle \zeta ,\zeta \rangle = \sum_{i,j} \langle \xi_i , \langle \eta_i ,\eta_j \rangle \xi_j \rangle = \langle \xi , X\xi \rangle$$
where $\xi = (\xi_1, \ldots , \xi_n )\in {\cal E}_{A_1}(B)\oplus \cdots \oplus {\cal E}_{A_n}(B)$ and $X$ is the matrix
$$X = \left( \begin{array}{ccc}
\langle \eta_1 ,\eta_1 \rangle & \cdots & \langle \eta_1 , \eta_n \rangle \\
\vdots & \ddots & \vdots \\
\langle \eta_n ,\eta_1 \rangle & \cdots & \langle \eta_n , \eta_n \rangle \\
\end{array} \right)$$
belonging to the $C^\ast$-algebra $\Hom (A_1\oplus \cdots \oplus A_n,A_1\oplus \cdots \oplus A_n)$.

By lemma \ref{pretens2} the matrix $X$ is positive.  By lemma \ref{pretens1} the product $\langle \zeta ,\zeta \rangle$ is therefore positive.
\end{proof}

For each object $B\in \Ob ({\cal B})$ form the vector space
$$N(B) = \{ \eta \in {\cal E}\odot_{\cal A}{\cal F}(B,-) \ |\ \langle \eta ,\eta \rangle = 0 \}$$
Then the collection of quotients ${\cal E}\odot_{\cal A}{\cal F}(B,-) /N(B)$ is a right $\cal B$-module equipped with an inner product.  It can be graded if desired.

\begin{definition}
We define the {\em inner tensor product}, ${\cal E}\gtp_{\cal A}{\cal F}$, to be the Hilbert $\cal B$-module obtained by completing the above inner product $\cal B$-module.
\end{definition}

If $T\colon {\cal E}\rightarrow {\cal E}'$ is a bounded operator then we define a bounded operator $T\otimes 1\colon {\cal E}\gtp_{\cal A}{\cal F}\rightarrow {\cal E}'\gtp_{\cal A}{\cal F}$ by the formula
$$(T\otimes 1) (\eta \otimes \xi) = (T\eta ) \otimes \xi$$

For example, let $F\colon {\cal A}\rightarrow {\cal B}$ be a $C^\ast$-functor, and let $\cal E$ be a Hilbert $\cal A$-module.  Then by example \ref{onedim} the category $\cal B$ is itself a Hilbert $({\cal A},{\cal B})$-bimodule and we can form an inner tensor product which we write ${\cal E}\gtp_F{\cal B}$.

The following fact is useful in computations.

\begin{proposition} \label{generators}
Consider a vector $\zeta \in {\cal E}\odot_{\cal A}{\cal F}(B,-)$.  Then the vector $\zeta$ belongs to the space $N(B)$ if and only if we can write
$$\zeta = \sum_{i=1}^n (\eta_i x_i \otimes \xi_i - \eta_i \otimes x_i \xi_i)$$
for vectors $\eta_i \in {\cal E}(A_i')$, $\xi_i \in {\cal F}(B,A_i)$, and morphisms $x_i\in \Hom (A_i,A_i')_{\cal A}$.
\end{proposition}

\begin{proof}
Let $\zeta = \eta x \otimes \xi - \eta \otimes x\xi$.  Then a straightforward calculation tells us that $\langle \zeta ,\zeta \rangle =0$.  Conversely, suppose that $\zeta \in N(B)$.  Write
$$\zeta = \eta_1\otimes \xi_1 + \cdots + \eta_n\otimes \xi_n$$
for vectors $\eta_i \in {\cal E}(A_i)$ and $\xi_i \in {\cal F}(B,A_i)$.  Form the matrix
$$X = \left( \begin{array}{ccc}
\langle \eta_1 ,\eta_1 \rangle & \cdots & \langle \eta_1 , \eta_n \rangle \\
\vdots & \ddots & \vdots \\
\langle \eta_n ,\eta_1 \rangle & \cdots & \langle \eta_n , \eta_n \rangle \\
\end{array} \right)$$
and write $\xi = (\xi_1 , \ldots ,\xi_n )\in {\cal F}(B,A_1)\oplus \cdots \oplus {\cal F}(B,A_n)$.  Then
$$\langle \xi , X\xi \rangle = \langle \zeta ,\zeta \rangle = 0$$

By lemma \ref{pretens2} the matrix $X$ is positive, so we can find a positive square root $X^\frac{1}{2}$.  We have the formula
$$\langle X^\frac{1}{2} \xi ,X^\frac{1}{2} \xi \rangle = \langle  \xi ,X\xi \rangle =0$$
and so $X^\frac{1}{2}\xi =0$.  Repeating this process we see that $X^\frac{1}{4} \xi =0$.
 
Form the vector $\eta = (\eta_1 ,\ldots ,\eta_n ) \in {\cal E}(A_1)\oplus \cdots \oplus {\cal E}(A_n)$.  Then by definition of the matrix $X$ we can write $X = \langle \eta ,\eta \rangle$.  By lemma \ref{pretens3} we can find a vector $\eta = (\eta_1', \ldots ,\eta_n')\in {\cal E}(A_1)\oplus \cdots \oplus {\cal E}(A_n)$ such that $\eta' X^\frac{1}{4} = \eta$.

Write
$$X^\frac{1}{4} = \left( \begin{array}{ccc}
x_{1,1} & \cdots & x_{1,n} \\
\vdots & \ddots & \vdots \\
x_{n,1} & \cdots & x_{n,n} \\
\end{array} \right)$$

Then we have the formulae
$$\sum_j x_{i,j} \xi_j = 0 \qquad \eta_j = \sum_i \eta_i' x_{i,j}$$
and so we conclude that
$$\zeta = \sum_{i,j} (\eta_i'x_{i,j}\otimes \xi_j - \eta_i'\otimes x_{i,j}\xi_j )$$
as required.
\end{proof}

Thus the Hilbert $\cal B$-module ${\cal E}\gtp_{\cal A}{\cal F}$ is generated by elementary tensors $\eta \otimes \xi$ satisfying the relation
$$\eta x \otimes \xi = \eta \otimes x \xi$$

Our final lemma is an application of the Kasparov stabilisation theorem.

\begin{lemma} \label{count}
Suppose that $\cal A$ and $\cal B$ are $\sigma$-unital graded $C^\ast$-categories.  Let $\cal E$ be a countably generated graded Hilbert $\cal A$-module and let $F\colon {\cal A}\rightarrow {\cal B}$ be a graded unital $C^\ast$-functor.  Then the inner tensor product ${\cal E}\gtp_F {\cal B}$ is countably generated.
\end{lemma}

\begin{proof}
Let $\cal H$ denote the Hilbert space of all sequences $(\alpha_n )$ in the field $\mathbb F$ such that the sum $\sum |\alpha_n |^2$ converges.  Then the Hilbert $\cal A$-module ${\cal H}(-,A)_{\cal A}$ is isomorphic to the outer tensor product ${\cal H}\gtp \Hom (-,A)_{\cal A}$.

Hence, by proposition \ref{generators}, the vector space ${\cal H}(-,A)_{\cal A}\gtp_F \Hom(B,-)_{\cal B}$ is generated by elements of the form
$$((\alpha_n)\otimes x)\otimes y = (\alpha_n)\otimes F(x)y$$
where $(\alpha_n ) \in {\cal H}$, $x\in \Hom (A',A)_{\cal A}$, and $y\in \Hom (B,F(A'))_{\cal B}$.  Since the $C^\ast$-functor $F$ is unital, every vector of the form $(\alpha_n)\otimes y$, where $y\in \Hom (B,F(A))_{\cal B}$, is obtained in this fashion.  Therefore we have isomorphisms
$${\cal H}(-,A)_{\cal A}\gtp_F {\cal B}\cong {\cal H}\gtp \Hom(-,F(A))_{\cal B}\cong {\cal H}(-,F(A))_{\cal B}$$

Let $\cal E$ be a countably generated Hilbert $\cal A$-module.  By the
Kasparov stabilisation theorem the Hilbert module $\cal E$ can be
considered to be a submodule of the
direct sum of countably many Hilbert $\cal A$-modules of the form ${\cal H}(-,A)_{\cal A}$.  Hence the above calculation tells us that the inner tensor product ${\cal E}\gtp_F {\cal B}$ is countably generated.
\end{proof}

\subsection{Some Computations}

We present in this subsection some computations concerning the
$C^\ast$-categories ${\cal L}({\cal A})$, ${\cal K}({\cal A})$, and
${\cal Q}({\cal A})$ which will help us to deduce some of the
properties of $KK$-theory.  

\begin{proposition}
The assignments ${\cal A}\mapsto {\cal L}({\cal A})$ and ${\cal A}\mapsto {\cal Q}({\cal A})$ are covariant functors from the category of $\sigma$-unital graded $C^\ast$-categories to the category of unital graded $C^\ast$-categories.\footnote{In fact there are set-theoretic difficulties involved here since for a graded $C^\ast$-category $\cal A$ the collection of all countably generated graded Hilbert $\cal A$-modules is not a set, even when the $C^\ast$-category $\cal A$ is itself small.  We can evade this problem by insisting that all of our Hilbert modules over a particular $C^\ast$-category lie within a given universe.}
\end{proposition}

\begin{proof}
Let $F\colon {\cal A}\rightarrow {\cal B}$ be a graded $C^\ast$-functor.  Let $\cal E$ be a countably generated graded Hilbert $\cal A$-module.  Then by lemma \ref{count} we have a countably generated graded Hilbert $\cal B$-module defined by the formula $F_\star ({\cal E}) = {\cal E}\gtp_F {\cal B}$.

If $T\colon {\cal E}\rightarrow {\cal F}$ is an operator between countably generated Hilbert $\cal B$-modules there is an induced operator $T\otimes 1 \colon {\cal E}\gtp_F {\cal B}\rightarrow {\cal F}\gtp_F {\cal B}$.  The operator $T\otimes 1$ is compact if the operator $T$ is compact.

It is now easy to check that the above process defines a $C^\ast$-functor $F_\star \colon {\cal L}({\cal A})\rightarrow {\cal L}({\cal B})$, that $1_\star =1$, and that $(FG)_\star = F_\star G_\star$.  The assignments ${\cal A}\mapsto {\cal L}({\cal A})$ and ${\cal A}\mapsto {\cal Q}({\cal A})$ are therefore covariant functors as claimed.
\end{proof}

\begin{definition}
We write ${\cal L}({\cal H}_{\cal A})$ to denote the full subcategory of the $C^\ast$-category ${\cal L}({\cal A})$ in which the objects are Hilbert $\cal A$-modules of the form ${\cal H}(-,A)_{\cal A}$ where $A\in \Ob ({\cal A})$.  We write ${\cal K}({\cal H}_{\cal A})$ to denote the corresponding full subcategory of the $C^\ast$-category ${\cal K}({\cal A})$.
\end{definition}

By proposition \ref{kideall} the $C^\ast$-category ${\cal K}({\cal H}_{\cal A})$ is a $C^\ast$-ideal in the $C^\ast$-category ${\cal L}({\cal H}_{\cal A})$.  We may therefore form the quotient ${\cal Q}({\cal H}_{\cal A}) = {\cal L}({\cal H}_{\cal A}) / {\cal K}({\cal H}_{\cal A})$

\begin{proposition}
Let $\cal A$ be a unital $C^\ast$-category which is equivalent to some $C^\ast$-algebra.  Then we have equivalences of $K$-theory spectra
$${\mathbb K}{\cal L}({\cal A})\simeq {\mathbb K}{\cal L}({\cal H}_{\cal A}) \qquad {\mathbb K}{\cal K}({\cal A})\simeq {\mathbb K}{\cal K}({\cal H}_{\cal A}) \qquad {\mathbb K}{\cal Q}({\cal A})\simeq {\mathbb K}{\cal Q}({\cal H}_{\cal A})$$
\end{proposition}

\begin{proof}
This result is an immediate consequence of the Kasparov stabilisation theorem and the definition of the $K$-theory groups.
\end{proof}

\begin{proposition} \label{swine}
The $K$-theory groups of the graded $C^\ast$-category ${\cal L}({\cal A})$ are all trivial.
\end{proposition}

\begin{proof}
We will use an `Eilenberg swindle' argument to prove this result.  

Let $\mathbb F_{p,q}$ be any Clifford algebra.  Let $\cal E$ be a
countably generated graded Hilbert $\cal A$-module, and let $x\in
SS({\cal L}({\cal E})\gtp {\mathbb F}_{p,q})$ be a supersymmetry.

Consider the direct sum, ${\cal E}^\infty$, of (countably) infinitely
many copies of the Hilbert $\cal A$-module $\cal E$.  We can form a supersymmetry
$$x\oplus x\oplus \cdots \in SS({\cal L}({\cal E}^\infty )\gtp {\mathbb F}_{p,q})$$

At the level of $K$-theory
$$\langle x\rangle + \langle x\oplus x\oplus \cdots \rangle = \langle x\oplus x\oplus \cdots \rangle$$
and so $\langle x \rangle = 0$.  Hence $K_1 ({\cal L}({\cal A})\gtp {\mathbb F}_{p,q}) = 0$ for all natural numbers $p$ and $q$.  By the Bott periodicity theorem we deduce that $K_n {\cal L}({\cal A})=0$ for all integers $n\in \mathbb Z$.
\end{proof}

\begin{proposition} \label{awkward}
The $C^\ast$-categories ${\cal B}\gtp {\cal K}$ and ${\cal K}({\cal
  H}_{\cal B})$ are isomorphic.
\end{proposition}

\begin{proof}
Let ${\cal K}({\cal B}^n )$ be the $C^\ast$-category of compact
operators between Hilbert $\cal B$-modules of the form $\Hom (-,
A)_{\cal B}^n$.  Then the $C^\ast$-category ${\cal K}({\cal H}_{\cal
  B})$ is the direct limit of the $C^\ast$-categories ${\cal K}({\cal
  B}^n)$.  Similarly, the $C^\ast$-category ${\cal B}\gtp {\cal K}$
is the direct limit of the matrix $C^\ast$-categories ${\cal B}\gtp
M_n ({\mathbb F})$.  It therefore suffices to show that the
$C^\ast$-categories ${\cal B}\gtp M_n ({\mathbb F})$ and ${\cal
  K}({\cal B}^n )$ are isomorphic.

We can define a graded $C^\ast$-functor $F\colon {\cal K}({\cal B}^n
)\rightarrow {\cal B}\gtp M_n ({\mathbb F})$ by mapping the Hilbert
$\cal B$-module $\Hom (-,A)_{\cal B}^n$ to the object $A$, and the
operator
$$\eta \mapsto (x_1 \oplus \cdots \oplus x_n ) \langle y_1 \oplus
\cdots \oplus y_n , \eta \rangle$$
to the matrix
$$\left( \begin{array}{ccc}
x_1 y_1^\star & \cdots & x_1 y_n^\star \\
\vdots & \ddots & \vdots \\
x_n y_1^\star & \cdots & x_n y_n^\star \\
\end{array} \right) \in \Hom (B,C)_{{\cal B}\gtp M_n ({\mathbb
  F})}$$
where $x_i \in \Hom (A,C)$ and $y_i \in \Hom (A,B)$.

Let $\{ e_\lambda \ |\ \lambda \in \Lambda \}$ be a self-adjoint
approximate unit for the $C^\ast$-algebra $\Hom (B,B)_{\cal
  B}$.
A simple calculation shows that
$$\inf_{\lambda \in \Lambda} \{ \| xe_\lambda - x \| \ |\ \lambda
\in \Lambda \} = 0$$
for any morphism $x\in \Hom (B,C)_{\cal B}$.

It follows that the $C^\ast$-functor $F$ has dense range.  As we
mentioned in section \ref{cstarcat}, the image of any $C^\ast$-functor
is closed.  Therefore the $C^\ast$-functor $F$ is surjective.

Now, consider the compact operator
$$T\colon \eta \mapsto \sum_i (x_1^{(i)} \oplus \cdots \oplus x_n^{(i)} ) \langle y_1^{(i)} \oplus
\cdots \oplus y_n^{(i)} , \eta \rangle$$
where $x_j^{(i)} \in \Hom (A,C)$ and $y_j^{(i)} \in \Hom (A,B)$.  Let
$\eta = \eta_1 \oplus \cdots \oplus \eta_n$.  Then
$$T\eta = \sum_{i,j} x_1^{(i)} {y_j^{(i)}}^\star \eta_j \oplus
\cdots \oplus x_n^{(i)} {y_j^{(i)}}^\star \eta_j$$

Thus
$$\| T \| = \sup_{\| \eta \|\leq 1} \left\| \sum_{i,j} x_1^{(i)} {y_j^{(i)}}^\star \eta_j \oplus
\cdots \oplus x_n^{(i)} {y_j^{(i)}}^\star \eta_j \right\|$$

By taking the various vectors $\eta_j$ to be elements of approximate
units, we see that the norm, $\| T\|$, is equal to the norm of the
matrix
$$\left( \begin{array}{ccc}
x_1 y_1^\star & \cdots & x_1 y_n^\star \\
\vdots & \ddots & \vdots \\
x_n y_1^\star & \cdots & x_n y_n^\star \\
\end{array} \right)$$

Hence the $C^\ast$-functor $F$ is an isometry, and we are done.
\end{proof}

\begin{theorem} \label{qee}
Let $\cal A$ be a small unital graded $C^\ast$-category which is equivalent to some $C^\ast$-algebra.  Then we have a natural stable equivalence of spectra
$$\Omega {\mathbb K}{\cal Q}({\cal A})\simeq {\mathbb K}({\cal A})$$
\end{theorem}

\begin{proof}
We have a natural short exact sequence
$$0\rightarrow {\cal K}({\cal A})\rightarrow {\cal L}({\cal A})\rightarrow {\cal Q}({\cal A})\rightarrow 0 $$
and so by proposition 4.7 of \cite{Mitch2.5} a natural fibration
$${\mathbb K}{\cal K}({\cal A})\rightarrow {\mathbb K}{\cal L}({\cal A})\rightarrow {\mathbb K}{\cal Q}({\cal A})$$

We therefore have an induced long exact sequence
$$\rightarrow K_{n+1}{\cal L}({\cal A})\rightarrow K_{n+1}{\cal Q}({\cal A})\stackrel{{\partial}_\star}{\rightarrow} K_n{\cal K}({\cal A})\rightarrow K_n{\cal L}({\cal A})\rightarrow$$
where the boundary map $\partial_\star \colon K_{n+1}{\cal Q}({\cal A})\rightarrow K_n{\cal K}({\cal A})$ is induced from a natural map $\partial \colon \Omega {\mathbb K}{\cal Q}({\cal A})\rightarrow {\mathbb K}{\cal K}({\cal A})$ defined by the homotopy lifting property of a fibration.

By proposition \ref{swine} the $K$-theory groups $K_{n+1}{\cal L}({\cal A})$ and $K_n{\cal L}({\cal A})$ are both zero.  Hence the map $\partial \colon \Omega {\mathbb K}{\cal Q}({\cal A})\rightarrow {\mathbb K}{\cal K}({\cal A})$ is a stable equivalence of spectra.

But by theorem \ref{bigstab} the inclusion ${\cal A}\hookrightarrow {\cal A}\gtp {\cal K}$ induces isomorphisms of $K$-theory groups, so by proposition \ref{awkward} there is a natural stable equivalence of spectra ${\mathbb K}{\cal K}({\cal H}_{\cal A})\simeq {\mathbb K}({\cal A})$.  It follows that there is a natural stable equivalence of spectra $\Omega {\mathbb K}{\cal Q}({\cal A})\simeq {\mathbb K}({\cal A})$.
\end{proof}

\section{$KK$-theory}

\subsection{Definition of $KK$-theory}

We are now ready to generalise the definitions and some results
concerning the $KK$-theory of $C^\ast$-algebras to the world of small
$C^\ast$-categories.  For many results concerning $KK$-theory to be
valid we need to restrict our attention to $\sigma$-unital
$C^\ast$-categories.  We will mention explicitly when the assumption
of $\sigma$-unitality is needed.

\begin{definition}
Let $\cal A$ and $\cal B$ be small graded $C^\ast$-categories. Then a {\em Kasparov $({\cal A},{\cal B})$-cycle} consists of a Hilbert $({\cal A},{\cal B})$-bimodule $\cal E$ together with a collection of degree $1$ operators $F_A\colon {\cal E}(-,A)\rightarrow {\cal E}(-,A)$ such that the norm
$$\| F \| = \sup \{ \| F_A  \| \ |\ A\in \Ob ({\cal A}) \}$$
is finite and the operators
$$x(F_A-F_A^\star) \qquad x(F_A^2 - 1) \qquad xF_A - (-1)^{\deg (x)} F_Bx$$
are compact for all morphisms $x\in \Hom (A,B)_{\cal A}$.
\end{definition}

An element of $KK$-theory is defined to be a certain equivalence class of Kasparov cycles.  We will abuse notation slightly, and write $F$ to denote both the collection of operators $F_A\colon {\cal E}(-,A)\rightarrow {\cal E}(-,A)$ and an individual operator of the form $F_A$.

\begin{definition}
Let $({\cal E},F)$ and $({\cal E}',F')$ be Kasparov $({\cal A},{\cal B})$-cycles.  Then the {\em direct sum} is the Kasparov cycle
$$({\cal E},F)\oplus ({\cal E}',F') = ({\cal E}\oplus {\cal E}',F\oplus F')$$
\end{definition}

\begin{definition}
\mbox{}

\begin{itemize}

\item[$\bullet$] A Kasparov $({\cal A},{\cal B})$-cycle $({\cal E},F)$ is called {\em degenerate} if the operators
$$x(F-F^\star) \qquad x(F^2 - 1) \qquad xF - (-1)^{\deg (x)} Fx$$
are equal to zero for all morphisms $x\in \Hom (A,B)_{\cal A}$.

\item[$\bullet$] An {\em operator homotopy} between Kasparov $({\cal A},{\cal B})$-cycles $({\cal E},F)$ and $({\cal E},F')$ is a norm-continuous path $({\cal E},F_t)$ of Kasparov $({\cal A},{\cal B})$-cycles such that $F_0 = F$ and $F_1 = F'$.

\item[$\bullet$] Kasparov $({\cal A},{\cal B})$-cycles $({\cal E}_1,F_1)$ and $({\cal E}_2, F_2)$ are called {\em equivalent} if there are degenerate Kasparov $({\cal A},{\cal B})$-cycles $({\cal E}_1',F_1')$ and $({\cal E}_2',F_2')$ such that the direct sums $({\cal E}_1,F_1)\oplus({\cal E}_1',F_1')$ and $({\cal E}_2, F_2)\oplus({\cal E}_2', F_2')$ are operator homotopic.

\end{itemize}

\end{definition}

We write $[({\cal E}, F)]$ to denote the equivalence class of a Kasparov $({\cal A},{\cal B})$-cycle $({\cal E},F)$, and $KK({\cal A},{\cal B})$ to denote the set of equivalence classes.

\begin{theorem}
The set $KK({\cal A},{\cal B})$ is an Abelian group with an operation defined by taking the direct sum of Kasparov cycles.  
\end{theorem}

\begin{proof}
It is easy to check that the set $KK({\cal A},{\cal B})$ is an Abelian semigroup, with identity element $[({\cal E},F)]$ where $({\cal E},F)$ is any degenerate Kasparov $({\cal A},{\cal B})$-cycle.

If ${\cal E}$ is a graded Hilbert $\cal B$-module, with grading ${\cal E}(A)= {\cal E}(A)_0 \oplus {\cal E}(A)_1$, define ${\cal E}^\mathrm{op}$ to be the Hilbert $\cal B$-module with the opposite grading.  Given a Kasparov $({\cal A},{\cal B})$ cycle $({\cal E}, F)$ define $(\tilde{\cal E}, -F)$ to be the Kasparov cycle in which the bimodule $\tilde{\cal E}$ is defined by writing $\tilde{\cal E}(-,A) = {\cal E}(-,A)^\mathrm{op}$ for all objects $A\in \Ob ({\cal A})$ and $\tilde{\cal E}(x_0 + x_1 ) = {\cal E}(x_0 -x_1)$ for all morphisms $x_0 , x_1 \in \Hom (A,A')_{\cal A}$ of degrees $0$ and $1$ respectively.

We can define an operator homotopy between the Kasparov cycle $({\cal E},F)\oplus (\tilde{\cal E},-F)$ and the cycle $\left( {\cal E}\oplus \tilde{\cal E} , \left( \begin{array}{cc}
0 & 1 \\
1 & 0 \\
\end{array} \right) \right)$ by the formula
$$G_\theta = \left( \begin{array}{cc}
F\cos \theta & \sin \theta \\
\sin \theta & -F\cos \theta \\
\end{array} \right) \qquad \theta \in [0, \frac{\pi}{2}]$$

But the cycle $\left( {\cal E}\oplus \tilde{\cal E} , \left( \begin{array}{cc}
0 & 1 \\
1 & 0 \\
\end{array} \right) \right)$ is degenerate so we have proved that
$$[({\cal E},F)] + [(\tilde{\cal E}, -F)] = 0$$

Therefore the set $KK({\cal A},{\cal B})$ is an abelian group.
\end{proof}

It is easy to see that when $A$ and $B$ are graded $C^\ast$-algebras the group $KK(A,B)$ defined above is the usual $KK$-theory group associated to graded $C^\ast$-algebras.

\begin{proposition} \label{kkfunc}
Let $\cal A$ and $\cal B$ be small $\sigma$-unital graded
$C^\ast$-categories.  Then the abelian group $KK({\cal A},{\cal B})$ is contravariantly functorial in the variable $\cal A$ and covariantly functorial in the variable $\cal B$.
\end{proposition}

\begin{proof}
Let $({\cal E},F)$ be a Kasparov $({\cal A},{\cal B})$-cycle and let
$G\colon {\cal B}\rightarrow {\cal B}'$ be a $C^\ast$-functor.  Then
we can define a Hilbert $({\cal A},{\cal B})$-bimodule ${\cal
  E}\gtp_G{\cal B}'$ and a collection of operators $F_A\otimes 1\colon
{\cal E}\gtp_G{\cal B}'\rightarrow {\cal E}\gtp_F{\cal B}'$.  It is
straightforward to check that we have defined a Kasparov $({\cal
  A},{\cal B}')$-cycle\footnote{Each Hilbert $\cal B$-module ${\cal
    E}_A\gtp_G{\cal B}'$ is countably generated by lemma \ref{count}.
  Here we need the assumption that the $C^\ast$-category ${\cal B}'$
  is $\sigma$-unital.}
$$G_\star ({\cal E}, F) = ({\cal E}\gtp_G{\cal B}' , F\otimes 1)$$
and that with such induced cycles we have a functorially induced map $G_\star \colon KK({\cal A},{\cal B})\rightarrow KK({\cal A},{\cal B}')$.

Now suppose we have a $C^\ast$-functor $H\colon {\cal A}'\rightarrow {\cal A}$.  Then the composite $C^\ast$-functor ${\cal E}\circ H\colon {\cal A}'\rightarrow {\cal L}({\cal B})$ is a Hilbert $({\cal A'},{\cal B})$-bimodule.  We have a functorially induced map of $KK$-theory groups defined by the formula
$$H^\star ({\cal E},F) = ({\cal E}\circ H ,F)$$
\end{proof}

Note that the above proof does not work if we allow the $C^\ast$-category $\cal
B$ to be non-$\sigma$-unital.  

\subsection{Duality}
In \cite{Hi2} the $K$-homology of a $C^\ast$-algebra $A$ is defined in terms of the ordinary $K$-theory of a `dual algebra' constructed from $A$.  We can in fact extend this approach to define the groups $KK^{-n} ({\cal A},{\cal B})$ for graded $C^\ast$-categories $\cal A$ and $\cal B$ in terms of the ordinary $K$-theory of some `dual category'.  Having done this we will be able to deduce information concerning $KK$-theory from information concerning $K$-theory.

\begin{definition}
An {\em operator}, $T$, between graded Hilbert $({\cal A},{\cal B})$-bimodules $\cal E$ and $\cal F$ is a collection of bounded operators $T_A \colon {\cal E}(-,A) \rightarrow {\cal F}(-,A)$ such that the norm
$$\| T \| = \sup \{ \| T_A \| \ |\ A\in \Ob ({\cal A}) \}$$
exists.  The operator $T$ is said to be of degree $1$ if each operator $T_A$ is of degree $1$, and of degree $0$ if each operator $T_A$ is of degree $0$.
\end{definition}

Note that we do {\em not} require an operator between Hilbert $({\cal A},{\cal B})$-bimodules to be a natural transformation.

\begin{definition}
Let $\cal A$ and $\cal B$ be $\sigma$-unital small graded
$C^\ast$-categories.  We define the {\em dual $C^\ast$-category}, ${\cal D}({\cal A},{\cal B})$, of the categories $\cal A$ and $\cal B$ to be the $C^\ast$-category in which the objects are graded $({\cal A},{\cal B})$-bimodules and the morphisms between $({\cal A},{\cal B})$-bimodules $\cal E$ and $\cal F$ are operators $T$ such that the graded commutator $xT_A - (-1)^{\deg (x) \deg (T)}T_B x$ is compact for all morphisms $x\in \Hom (A,B)_{\cal A}$

We define the subcategory of {\em locally compact operators}, ${\cal KD}({\cal A},{\cal B})$, to be the subcategory of the category ${\cal D}({\cal A},{\cal B})$ in which the morphisms are operators $T\colon {\cal E}\rightarrow {\cal F}$ such that the operators $T_Bx$ and $xT_A$ are compact for all morphisms $x\in \Hom (A,B)_{\cal A}$.
\end{definition}

It is easy to check that the collection of objects and morphisms ${\cal D}({\cal A},{\cal B})$ is a unital graded $C^\ast$-category with direct sum, and that the collection ${\cal K}{\cal D}({\cal A},{\cal B})$ is a $C^\ast$-ideal.  We can form the quotient ${\cal QD}({\cal A},{\cal B}) = {{\cal D}({\cal A},{\cal B})}/{{\cal KD}({\cal A},{\cal B})}$.

The quotient ${\cal QD}({\cal A},{\cal B})$ is a graded unital $C^\ast$-category.  The following result is proved similarly to proposition \ref{kkfunc}

\begin{proposition} 
The $C^\ast$-category ${\cal QD}({\cal A},{\cal B})$ is covariantly functorial in the variable $\cal B$, and covariantly functorial in the variable $\cal A$.
\noproof
\end{proposition}

The above proposition is not true if we allow the $C^\ast$-category
$\cal B$ to be non-$\sigma$-unital.

\begin{theorem} \label{dual}
There is a natural isomorphism $K_1 {\cal QD}({\cal A},{\cal B})\cong KK({\cal A},{\cal B})$.
\end{theorem}

\begin{proof}
The additive completion ${\cal QD}({\cal A},{\cal B})_\oplus$ is equivalent to the $C^\ast$-category ${\cal QD}({\cal A},{\cal B})$ in the sense of definition \ref{natiso2}, so the group $K_1 {\cal QD}({\cal A},{\cal B})$ is defined to be the set of formal differences
$$\{ \langle x\rangle  - \langle y\rangle \ |\ x,y\in SS({\cal E}),\ {\cal E}\in \Ob ({\cal QD}({\cal A},{\cal B})) \}$$

It is easy to check that we can define a natural homomorphism $\alpha \colon K_1{\cal QD}({\cal A},{\cal B})\rightarrow KK({\cal A},{\cal B})$ by writing
$$\alpha (\langle x\rangle -\langle y\rangle ) = [({\cal E},x)] - [({\cal E},y)]$$
for supersymmetries $x,y\in SS({\cal E})$.

Let $({\cal E},F)$ be any Kasparov cycle.  Then we can find degenerate cycles $({\cal E}',F')$ and $({\cal E}\oplus {\cal E}',G)$ for some Hilbert $({\cal A},{\cal B})$-bimodule ${\cal E}'$.  Hence
$$[({\cal E},F)] = [({\cal E}\oplus {\cal E}', F\oplus F')] - [({\cal E}\oplus {\cal E}',G)] = \alpha (\langle F\oplus F' \rangle - \langle G \rangle)$$
and we have proved that the homomorphism $\alpha$ is surjective.

Now suppose that $({\cal E},F)$ and $({\cal E},F')$ are Kasparov cycles and that $[({\cal E},F)] = [({\cal E},F')]$.  Then we can find a degenerate cycle $({\cal F},G)$ such that the cycles $({\cal E}\oplus {\cal F}, F\oplus G)$ and $({\cal E}\oplus {\cal F},F'\oplus G)$ are operator-homotopic.  Hence the supersymmetries $F\oplus G$ and $F'\oplus G$ lie in the same path-component of the space $SS({\cal E}\oplus {\cal F})$.

Thus
$$\begin{array}{rrcl}
& \langle F\oplus G \rangle & = & \langle F'\oplus G \rangle \\
\Rightarrow & \langle F \rangle & = & \langle F' \rangle \\
\end{array}$$
and the homomorphism $\alpha$ is also injective.
\end{proof}

\begin{lemma}
Suppose that $\cal A$ and $\cal B$ are unital graded $C^\ast$-categories that are equivalent to $C^\ast$-algebras.  Let $p,q\in \mathbb N$.  Then there is a natural isomorphism
$$K_1{\cal QD}({\cal A},{\cal B}\gtp \mathbb F_{p,q})\cong K_1({\cal QD}({\cal A},{\cal B})\gtp \mathbb F_{p,q})$$
\end{lemma}

\begin{proof}
Consider the Clifford algebra $\mathbb F_{1,0}$ generated by an element $e$ such that $e^2 =1$.  We can define a natural $C^\ast$-functor $\alpha \colon {\cal QD}({\cal A},{\cal B})\gtp \mathbb F_{1,0}\rightarrow {\cal QD}({\cal A},{\cal B}\gtp \mathbb F_{1,0})$ by mapping the Hilbert $({\cal A},{\cal B})$-bimodule $\cal E$ to the Hilbert $({\cal A},{\cal B}\gtp \mathbb F_{1,0})$-bimodule ${\cal E}\gtp \mathbb F_{1,0}$ and by mapping the morphism $x\otimes 1 + y\otimes e$, where $x,y\in \Hom ({\cal E},{\cal E}')_{{\cal QD}({\cal A},{\cal B})}$, to the operator $T$ defined by the formula
$$T_A(\eta \otimes 1 + \xi \otimes e) = x(\eta )\otimes 1 + x(\xi )\otimes e + (-1)^{\deg (\eta)}y(\eta )\otimes e + (-1)^{\deg (\xi )} y(\xi )\otimes 1$$

It is clear that the $C^\ast$-functor $\alpha$ is faithful, with image the full subcategory of the $C^\ast$-category ${\cal QD}({\cal A},{\cal B}\gtp \mathbb F_{1,0})$ in which the set of objects is the set of Hilbert $({\cal A},{\cal B})$-bimodules of the form ${\cal E}\gtp \mathbb F_{1,0}$, where $\cal E$ is a Hilbert $({\cal A},{\cal B})$-bimodule.  By the Kasparov stabilisation theorem we obtain a natural isomorphism
$$K_1{\cal QD}({\cal A},{\cal B}\gtp \mathbb F_{1.0})\cong K_1({\cal QD}({\cal A},{\cal B})\gtp \mathbb F_{1,0})$$

A similar calculation yields a natural isomorphism
$$K_1{\cal QD}({\cal A},{\cal B}\gtp \mathbb F_{0,1})\cong K_1({\cal QD}({\cal A},{\cal B})\gtp \mathbb F_{0,1})$$

The desired result now follows by induction since we have isomorphisms ${\mathbb F}_{p,q}\gtp {\mathbb F}_{r,s}\cong {\mathbb F}_{p+r,q+s}$.
\end{proof}

Because of the above lemma, theorem \ref{dual}, and the Bott periodicity theorem it makes sense to define further $KK$-theory groups by the formula
$$KK^{p-q} ({\cal A},{\cal B}) = KK({\cal A},{\cal B}\gtp \mathbb F_{p,q})$$

We have natural isomorphisms
$$KK^{p-q}({\cal A},{\cal B})\cong K_1({\cal QD}({\cal A},{\cal B})\gtp \mathbb F_{p,q})\cong K_{1-(p-q)}{\cal QD}({\cal A},{\cal B})$$

We can use the above notion of `duality' to define a spectrum for $KK$-theory.  This idea comes from \cite{HPR}.

\begin{definition} \label{kkspec}
Let $\cal A$ and $\cal B$ be small graded $C^\ast$-categories.  Then we define the {\em $KK$-theory spectrum} of the $C^\ast$-categories $\cal A$ and $\cal B$ to be the spectrum
$${\mathbb K}{\mathbb K}({\cal A},{\cal B}) = \Omega {\mathbb K} {\cal QD}({\cal A},{\cal B})$$
\end{definition}

We can similarly define $K$-homology spectra.

\begin{definition}
Let $\cal A$ be a small graded $C^\ast$-category.  Then we define the {\em $K$-homology spectrum} of the $C^\ast$-category $\cal A$ to be the spectrum
$${\mathbb K}_\mathrm{hom}({\cal A}) = \Omega {\mathbb K} {\cal QD}({\cal A}, {\mathbb F})$$
\end{definition}

Suppose that $X$ is a locally compact Hausdorff space.  Then the above definition gives us a spectrum for locally-finite $K$-homology
$${\mathbb K}^\mathrm{lf}_\mathrm{hom}(X) = \Omega {\mathbb K} {\cal QD}(C_0(X), {\mathbb F})$$

A slight adjustment delivers a spectrum for ordinary rather than locally finite $K$-homology
$${\mathbb K}_\mathrm{hom}(X) = \lim_{\longrightarrow \atop K^\mathrm{compact}\subseteq X} {\mathbb K}^{\mathrm lf}_\mathrm{hom}(K)$$

Finally, let us note that we can recover the definition of $K$-theory
from that of $KK$-theory.

\begin{proposition}
Let $\cal B$ be a small $\sigma$-unital graded $C^\ast$-category.  Then we have a natural stable equivalence of spectra
$${\mathbb K}{\mathbb K}({\mathbb F},{\cal B})\simeq {\mathbb K}({\cal B})$$
\end{proposition}

\begin{proof}
Observe that
$${\mathbb K}{\mathbb K}({\mathbb F},{\cal B}) = \Omega {\mathbb
  K}{\cal Q}({\cal B})$$

The result now follows from theorem \ref{qee}.
\end{proof}

\subsection{Products}
The purpose of this section is to define a few special cases of the Kasparov product for the $KK$-theory of $C^\ast$-categories.  The products we consider in this section can be expressed as natural morphisms in the stable category of symmetric spectra.

We begin with an aside on trivially graded $C^\ast$-categories.

\begin{definition}
Let $\cal A$ be a small unital trivially graded $C^\ast$-category.  For each object $A\in \Ob ({\cal A}_\oplus)$ we form the set of projections
$$P(A) = \{ p\in \Hom (A,A) \ |\ p = p^\star, p^2 = p \}$$

We write $p\sim_h q$ when the projections $p\oplus 0_B\oplus 0_C$ and $0_A\oplus q\oplus 0_C$ lie in the same path-component of the space $P(A\oplus B\oplus C)$ for some object $C$.  We write $[p]$ for the equivalence class containing the projection $p$.
\end{definition}

Form the set $V_0({\cal A})$ of all equivalence classes of projections in the $C^\ast$-category ${\cal A}_\oplus$.  Then the set $V_0({\cal A})$ is an abelian semigroup with group operation induced by taking the direct sum, $p\oplus q$, of projections $p$ and $q$.  The following result comes from \cite{Mitch2.5}.  It is a consequence of the Bott periodicity theorem.

\begin{theorem} \label{isoku}
The $K$-theory group $K_0({\cal A})$ is naturally isomorphic to the Groth\-en\-dieck completion of the semigroup $V_0({\cal A})$.
\noproof
\end{theorem}

The first of the products we need to consider also comes from \cite{Mitch2.5}.

\begin{theorem} \label{product1}
There is a natural product ${\mathbb K}({\cal A})\wedge {\mathbb K}({\cal B})\rightarrow {\mathbb K}({\cal A}\gtp {\cal B})$ in the stable category of symmetric spectra.  When the $C^\ast$-categories $\cal A$ and $\cal B$ are trivially graded and unital, the induced map $K_0 ({\cal A})\otimes K_0({\cal B})\rightarrow K_0({\cal A}\otimes {\cal B})$ is defined by the formula
$$([p],[q])\mapsto [p\otimes q]$$
\noproof
\end{theorem}

The above product is called the {\em exterior product} of $K$-theory spectra.

\begin{proposition} \label{descent}
Let $\cal A$ and $\cal B$ be small $\sigma$-unital graded $C^\ast$-categories.  Then there is a canonical $C^\ast$-functor
$$P\colon {\cal A}\gtp {\cal QD}({\cal A},{\cal B})\rightarrow {\cal Q}({\cal B})$$

The $C^\ast$-functor $P$ is natural in the variable $\cal B$ in the obvious sense and natural in the variable $\cal A$ in the sense that if we have a $C^\ast$-functor $F\colon {\cal A}\rightarrow {\cal A}'$ then there is an induced commutative diagram
$$\begin{array}{ccc}
{\cal A}\gtp {\cal QD}({\cal A},{\cal B}) & \stackrel{P}{\longrightarrow} & {\cal Q}({\cal B}) \\
\uparrow & & \| \\
{\cal A}\gtp {\cal QD}({\cal A}',{\cal B}) & & {\cal Q}({\cal B}) \\
\downarrow & & \| \\
{\cal A}'\gtp {\cal QD}({\cal A}',{\cal B}) & \stackrel{P}{\longrightarrow} & {\cal Q}({\cal B}) \\
\end{array}$$
\end{proposition}

\begin{proof}
We can define such a $C^\ast$-functor by writing
$$P(A\otimes {\cal E}) = {\cal E}(-,A) \qquad P(x\otimes T) = xT_A$$
for all objects $A\otimes {\cal E}\in \Ob ({\cal A}\gtp {\cal QD}({\cal A},{\cal B}))$ and morphisms $x\otimes T\in \Hom (A\otimes {\cal E},B\otimes {\cal E}')_{{\cal A}\gtp {\cal QD}({\cal A},{\cal B})}$.
\end{proof}

\begin{corollary} \label{descent2}
There is a natural morphism
$${\mathbb K}({\cal A})\wedge {\mathbb K}{\mathbb K}({\cal A},{\cal B})\rightarrow {\mathbb K}({\cal B})$$
in the stable category of symmetric spectra
\end{corollary}

\begin{proof}
By the previous proposition and theorem \ref{qee} we have a natural morphism
$$\Omega {\mathbb K}({\cal A}\gtp {\cal QD}({\cal A},{\cal B}))\rightarrow \Omega{\mathbb K}({\cal Q}({\cal B}))\rightarrow {\mathbb K}({\cal B})$$

Composition with the exterior product gives us a natural morphism
$${\mathbb K}({\cal A})\wedge {\mathbb K}{\mathbb K}({\cal A},{\cal B})\rightarrow \Omega {\mathbb K}({\cal A}\gtp {\cal QD}({\cal A},{\cal B}))\rightarrow {\mathbb K}({\cal B})$$
as required.
\end{proof}

\begin{proposition} \label{diagonal}
Let $\cal A$, $\cal B$, and $\cal C$ be small $\sigma$-unital graded
$C^\ast$-categories.  Then there is a canonical $C^\ast$-functor
$$D\colon {\cal QD}({\cal A},{\cal B})\rightarrow {\cal QD}({\cal A}\gtp {\cal C},{\cal B}\gtp {\cal C})$$

The $C^\ast$-functor $D$ is natural in the variables $\cal A$ and $\cal B$ and natural in the variable $\cal C$ in the sense that for any $C^\ast$-functor $F\colon {\cal C}\rightarrow {\cal C}'$ there is an induced commutative diagram
$$\begin{array}{ccc}
{\cal QD}({\cal A},{\cal B}) & \stackrel{D}{\longrightarrow} & {\cal QD}({\cal A}\gtp {\cal C}, {\cal B}\gtp {\cal C}) \\
\| & & \downarrow \\
{\cal QD}({\cal A},{\cal B}) & & {\cal QD}({\cal A}\gtp {\cal C}, {\cal B}\gtp {\cal C}') \\
\| & & \uparrow \\
{\cal QD}({\cal A},{\cal B}) & \stackrel{D}{\longrightarrow} & {\cal QD}({\cal A}\gtp {\cal C}', {\cal B}\gtp {\cal C}') \\
\end{array}$$
\end{proposition}

\begin{proof}
We can define such a $C^\ast$-functor by writing
$$D({\cal E}) = {\cal E}\gtp {\cal C} \qquad D(T) = T\otimes 1$$
for all objects ${\cal E}\in \Ob ({\cal QD}({\cal A},{\cal B}))$ and morphisms $T\in \Hom ({\cal E},{\cal E}')_{{\cal QD}({\cal A},{\cal B})}$.
\end{proof}

We can use the above construction to define a `slant product' between $K$-theory and $K$-homology.

\begin{theorem} \label{slant}
There is a canonical morphism
$$S\colon {\mathbb K}({\cal A}\gtp {\cal B})\wedge {\mathbb K}_\mathrm{hom}({\cal A}) \rightarrow {\mathbb K}({\cal B})$$
in the stable category of symmetric spectra.  The morphism $S$ is natural in the variable $\cal B$ in the obvious sense and natural in the variable $\cal A$ in the sense that a $C^\ast$-functor $F\colon {\cal A}\rightarrow {\cal A}'$ induces a commutative diagram
$$\begin{array}{ccc}
{\mathbb K}({\cal A}\gtp {\cal B})\wedge {\mathbb K}_\mathrm{hom}({\cal A}) & \stackrel{S}{\longrightarrow} & {\mathbb K}({\cal B}) \\
\uparrow & & \| \\
{\mathbb K}({\cal A}\gtp {\cal B})\wedge {\mathbb K}_\mathrm{hom}({\cal A}') & & {\mathbb K}({\cal B}) \\
\downarrow & & \| \\
{\mathbb K}({\cal A}'\gtp {\cal B})\wedge {\mathbb K}_\mathrm{hom}({\cal A}') & \stackrel{S}{\longrightarrow} & {\mathbb K}({\cal B}) \\
\end{array}$$
\end{theorem}

\begin{proof}
As usual we write ${\mathbb K}_\mathrm{hom}({\cal A}) = {\mathbb K}{\mathbb K}({\cal A},\mathbb F)$.  By the above proposition and corollary \ref{descent2} we can form a morphism
$${\mathbb K}({\cal A}\gtp {\cal B})\wedge {\mathbb K}_\mathrm{hom}({\cal B})\rightarrow {\mathbb K}({\cal A}\gtp {\cal B})\wedge {\mathbb K}{\mathbb K}({\cal A}\gtp {\cal B},{\cal B})\rightarrow {\mathbb K}({\cal B})$$
\end{proof}

\section{Assembly}

\subsection{Finitely Generated Projective Hilbert Modules}

Let $A$ be a (trivially graded) unital $C^\ast$-algebra.  Recall that a Hilbert $A$-module, $\cal E$, is called {\em finitely generated and projective} if we can find a Hilbert $A$-module ${\cal E}'$ and an isomorphism of Hilbert modules $\phi \colon {\cal E}\oplus {\cal E}' \rightarrow A^n$ for some $n$.  A finitely generated Hilbert module defines an element of the $K$-theory group $K_0 (A)$.

We can generalise this idea to see how certain Hilbert modules over a (trivially graded) unital $C^\ast$-category $\cal A$ satisfying similar conditions defines a elements of the $K$-theory group $K_0 ({\cal A})$.

\begin{definition}
Let $\cal A$ be a small unital $C^\ast$-category.  Then a Hilbert $\cal
A$-module, $\cal E$, is called {\em finitely generated and projective}
if there is a Hilbert $\cal A$-module ${\cal E}'$ such that the direct sum ${\cal E}\oplus {\cal E}'$ is isomorphic to the direct sum of finitely many
Hilbert $\cal A$-modules of the form $\Hom (-,A)_{\cal A}$.
\end{definition}

We write ${\cal L}({\cal A}_\mathrm{fgp})$ to denote the category of finitely generated projective Hilbert $\cal A$-modules and bounded operators.

\begin{lemma}
The natural inclusion $i\colon {\cal A}\hookrightarrow {\cal L}({\cal A}_\mathrm{fgp})$ defined by mapping an object $A\in \Ob ({\cal A})$ to the Hilbert $\cal A$-module $\Hom (-,A)_{\cal A}$ induces a stable equivalence of $K$-theory spectra
$$i_\star \colon {\mathbb K}({\cal A})\rightarrow {\mathbb K}{\cal L}({\cal A}_\mathrm{fgp})$$
\end{lemma}

\begin{proof}
Consider supersymmetries $x,y\in SS({\cal E}\gtp {\mathbb F}_{p,q})$
for some finitely generated projective Hilbert $\cal A$-module $\cal
E$.  Since the module $\cal E$ is finitely generated and projective,
there is a Hilbert $\cal A$-module ${\cal E}'$ and objects $A_1 ,
\ldots , A_n \in \Ob ({\cal A})$ such that
$${\cal E}\oplus {\cal E}'\cong \Hom (-, A_1) \oplus \cdots \oplus
\Hom (-,A_n)$$

Choose a supersymmety $E\in SS({\cal E}'\gtp {\mathbb F}_{p,q})$.  Then, looking at $K$-theory classes:
$$\langle x\rangle - \langle y\rangle = \langle x\oplus E\rangle - \langle y\oplus E\rangle$$
which by the above isomorphism is the image of a $K$-theory class under the map induced by the $C^\ast$-functor $i$.  Hence the induced map
$$i_\star \colon K_1 ({\cal A}\gtp {\mathbb F}_{p,q})\rightarrow K_1 ({\cal L}({\cal A}_\mathrm{fgp})\gtp {\mathbb F}_{p,q})$$
is surjective.  Injectivity of the induced map $i_\star$ can be proved similarly.

By the Bott periodicity theorem it follows that the induced maps $i_\star \colon K_n ({\cal A})\rightarrow K_n {\cal L}({\cal A}_\mathrm{fgp})$ are all isomorphisms.  Hence the induced map of $K$-theory spectra
$$i_\star \colon {\mathbb K} ({\cal A})\rightarrow {\mathbb K} {\cal L}({\cal A}_\mathrm{fgp})$$
is a stable equivalence.
\end{proof}

Now, let $\cal A$ be a trivially graded unital $C^\ast$-category.   Then it is shown in \cite{Mitch2.5} that a projection, $p$, in the additive completion ${\cal A}_\oplus$ defines an element of the zeroth space in the $K$-theory spectrum:
$$[p]\in {\mathbb K}({\cal A})_0$$

If $\cal E$ is a finitely generated projective Hilbert $\cal A$-module, there is an associated projection $1_{\cal E}\colon {\cal E}\rightarrow {\cal E}$ in the category ${\cal L}({\cal A}_\mathrm{fgp})$.  By the above theorem we have an associated $K$-theory class
$$[{\cal E}] = i_\star^{-1} [1_{\cal E}] \in {\mathbb K}({\cal A})_0$$

By theorem \ref{slant} we have the following result.

\begin{proposition} \label{slant2}
Let $\cal A$ and $\cal B$ be trivially graded unital $C^\ast$-categories.  Let $\cal E$ is a finitely generated projective Hilbert ${\cal A}\otimes {\cal B}$-module.  Then there is an induced morphism of spectra:
$$[{\cal E}]\otimes \colon {\mathbb K}_\mathrm{hom}({\cal A}) \rightarrow {\mathbb K}({\cal B})$$

For all $C^\ast$-functors $F\colon {\cal A}\rightarrow {\cal A}'$ and $G\colon {\cal B}\rightarrow {\cal B}'$ there are commutative diagrams
$$\begin{array}{ccc}
{\mathbb K}_\mathrm{hom} ({\cal A}) & \stackrel{[{\cal E}]\otimes}{\longrightarrow} & {\mathbb K}({\cal B}) \\
\uparrow & & \| \\
{\mathbb K}_\mathrm{hom} ({\cal A}') & \stackrel{F_\star [{\cal E}]\otimes}{\longrightarrow} & {\mathbb K}({\cal B}) \\
\end{array}$$
and
$$\begin{array}{ccc}
{\mathbb K}_\mathrm{hom} ({\cal A}) & \stackrel{[{\cal E}]\otimes}{\longrightarrow} & {\mathbb K}({\cal B}) \\
\| & & \downarrow \\
{\mathbb K}_\mathrm{hom} ({\cal A}) & \stackrel{G_\star [{\cal E}]\otimes}{\longrightarrow} & {\mathbb K}({\cal B}') \\
\end{array}$$
\noproof
\end{proposition}

\subsection{The Analytic Assembly Map}

We are now ready to apply our machinery in order to give a natural description of the analytic assembly map at the level of spectra.  Such a description will enable us to use results from \cite{WW} to characterise the assembly map.

\begin{definition}
Let $X$ be a path-connected locally compact Hausdorff topological space that is locally path-connected.  Let $K\subseteq X$ be a compact subspace.  Choose a basepoint $x_0\in K$.  Then we form the space
$$(\tilde{K},x_0) = \frac{\{ \tilde{x} \in C([0,1]\rightarrow X) \ |\ \tilde{x} (0)\in K,\ \tilde{x} (1) = x_0 \} }{\sim}$$
where the equivalence relation $\sim$ is that of paths with fixed endpoints being homotopic in the space $X$.
\end{definition}

We write $p\colon (\tilde{K},x_0)\rightarrow X$ to denote the natural map defined by sending a path $\tilde{x}$ to the initial point $\tilde{x} (0)$.  The following result is easy to see.

\begin{proposition} \label{ucp}
Let $x\in K$.  Then there is an open subset $U\subseteq K$ such that $x\in U$ and we have a homeomorphism $p^{-1}[U]\approx U\times \pi_1 (X,x_0)$.
\noproof
\end{proposition}

Let $\tilde{x}_1,\tilde{x}_2 \colon [0,1]\rightarrow X$ be paths such that $\tilde{x}_1 (0) = \tilde{x}_2 (1)$.  Then we write $\tilde{x}_1 \tilde{x}_2$ to denote the path obtained by moving along the path $\tilde{x}_2$ followed by the path $\tilde{x}_1$.\footnote{This notation may appear `backwards', but is consistent with other conventions in this article such as writing $xy$ for the composition of a morphism $y$ with a morphism $x$ in some category.}  We can define an action of the fundamental group $\pi_1 (X,x_0 )$ on the space $(\tilde{K},x_0)$ by the formula $(g,\tilde{x} )\mapsto g\tilde{x}$.

Following \cite{Kas3} there is a canonical Hilbert $C(K)\otimes C^\star_\mathrm{max} \pi_1 (X,x_0)$-module based on this action.

\begin{proposition} \label{kasmod}
Define ${\cal F}_K (x_0)$ to be the space of all continuous maps $\mu \colon (\tilde{K},x_0)\rightarrow C^\star_\mathrm{max} \pi_1(X;x_0)$ such that $\mu (g\tilde{x} ) = g\mu (\tilde{x} )$ for all group elements $g\in \pi_1 (X;x_0 )$ and paths $\tilde{x}\in (\tilde{K},x_0)$.  Then the space ${\cal F}_K(x_0)$ is a finitely generated projective Hilbert $C(K)\otimes C^\star_\mathrm{max}\pi_1(X,x_0)$-module with inner product
$$\langle \mu , \mu' \rangle (p(\tilde{x} )) = \mu^\star (\tilde{x} ) \mu' (\tilde{x} )$$
and $C(K)\otimes C^\star_\mathrm{max}\pi_1 (X,x_0)$-action
$$\mu (f\otimes g)(\tilde{x} ) = f(p(\tilde{x} ))\mu (\tilde{x} ) g$$
\noproof
\end{proposition}

By theorem \ref{slant2} we can define a mapping
$$[{\cal F}_K (x_0) ]\otimes \colon K_\star (K)\rightarrow K_\star (C^\star_\mathrm{max} \pi_1 (X,x_0))$$

Recall that we can define the $K$-homology groups of the space $X$ to be the direct limit of the $K$-homology groups $K_\star (K)$ where $K$ is a compact subspace of $X$.  We therefore have a map
$$\beta_\star \colon K_\star (X)\rightarrow K_\star C^\star_\mathrm{max} \pi_1 (X,x_0)$$

\begin{definition}
The map $\beta_\star$ is called the {\em analytic assembly map}.
\end{definition}

We can generalise the above definition to form a basepoint-free map at the level of spectra.  The key idea is the following construction from \cite{Mitch2}.

\begin{proposition}
There is a functor ${\cal G}\mapsto C^\star_\mathrm{max} {\cal G}$ from the category of discrete groupoids to the category of unital $C^\ast$-categories such that

\begin{itemize}

\item[$\bullet$] There is a natural inclusion ${\cal G}\hookrightarrow C^\star_\mathrm{max}{\cal G}$

\item[$\bullet$] Let $f\colon {\cal G}\rightarrow {\cal H}$ be an equivalence of groupoids.  Then the induced map $f_\star \colon C^\star_\mathrm{max}{\cal G}\rightarrow C^\star_\mathrm{max}{\cal H}$ is an equivalence of $C^\ast$-categories

\item[$\bullet$] If $G$ is a group then the $C^\ast$-category $C^\star_\mathrm{max}G$ is the maximal $C^\ast$-algebra of the group $G$

\end{itemize}

\end{proposition}

The groupoid $C^\star_\mathrm{max}({\cal G})$ is called the {\em maximal $C^\ast$-category} of the groupoid $\cal G$.  It is constructed similarly to the reduced $C^\ast$-category of a groupoid that is considered in \cite{DL}.

Now, let $\pi (X)|_K$ be the full subgroupoid of the fundamental groupoid $\pi (X)$ in which the objects are the points of the subspace $K$.

Choose a point $x_0\in K$.  For each point $x\in K$ choose a path, $\gamma_x$, from the point $x_0$ to the point $x$.  Then there is a unitary element $u_x = 1\otimes [\gamma_x] \in \Hom (x_0 , x)_{C(K)\otimes C^\star_\mathrm{max}(\pi (X)|_K)}$ for each point $x\in K$.  The unital $C^\ast$-category $C(K)\otimes C^\star_\mathrm{max} \pi (X)|_K$ is therefore equivalent to the $C^\ast$-algebra $C(K)\otimes C^\star_\mathrm{max}(X,x_0)$.  Similarly the $C^\ast$-category $C^\star_\mathrm{max} \pi (X)|_K$ is equivalent to the $C^\ast$-algebra $C^\star_\mathrm{max} \pi (X,x_0)$.

\begin{definition}
For a point $x\in K$, define ${\cal F}_K(x)$ to be the set of maps
$$\{ \mu_y \colon (\tilde{K},y)\rightarrow \Hom (x,y)_{C^\star_\mathrm{max} \pi (X)|_K} \ |\ y\in K \}$$
such that
$$\mu_{y'}(g\tilde{y}) = g\mu_y (\tilde{y})$$
for all groupoid elements $g\in \Hom(y,y')_{\pi (X)|_K}$ and paths $\tilde{y}\in (\tilde{K},y)$.
\end{definition}

The collection, ${\cal F}_K$, of spaces of the form ${\cal F}_K(x)$ is a Hilbert $C(K)\otimes C^\star_\mathrm{max}\pi (X)|_K$-module with inner product
$$\langle \{ \mu_y \} , \{ \mu_y'\} \rangle (p(\tilde{y} )) = \mu_y^\star (\tilde{y} ) \mu_y' (\tilde{y} )$$
and $C(K)\otimes C^\star_\mathrm{max}\pi (X)|_K$-action
$$\mu (f\otimes g)(\tilde{y} ) = f(p(\tilde{y} ))\mu (\tilde{y} ) g$$

By proposition \ref{slant2} we obtain a morphism
$${\cal F}_K\otimes \colon {\mathbb K}_\mathrm{hom}(K)\rightarrow {\mathbb K}(C^\star_\mathrm{max} \pi (X)|_K)$$

At the level of homotopy groups we have an induced map
$$[{\cal F}_K]\otimes \colon K^\star (K)\rightarrow K_\star (C^\star_\mathrm{max}\pi (X)|_K)$$

\begin{lemma} \label{biglem1}
We have a commutative diagram
$$\begin{array}{ccc}
K_\star (K) & \stackrel{[{\cal F}_K (x_0)]\otimes}{\longrightarrow} & K_\star C^\star_\mathrm{max} \pi_1 (X,x_0) \\
\| & & \downarrow \\
K_\star (K) & \stackrel{[{\cal F}]\otimes}{\longrightarrow} & K_\star C^\star_\mathrm{max} \pi (X)|_K \\
\end{array}$$
such that the vertical map on the right is an isomorphism.
\end{lemma}

\begin{proof}
The result follows immediately from the fact that the $C^\ast$-category $C(K)\otimes C^\star_\mathrm{max} \pi (X)|_K$ is equivalent to the $C^\ast$-algebra $C(K)\otimes C^\star_\mathrm{max}(X,x_0)$ and the $C^\ast$-category $C^\star_\mathrm{max} \pi (X)|_K$ is equivalent to the $C^\ast$-algebra $C^\star_\mathrm{max} \pi (X,x_0)$.
\end{proof}

\begin{lemma} \label{biglem2}
Let $s\colon (X,K)\rightarrow (Y,L)$ be a map of pairs of spaces, where the spaces $X$ and $Y$ are path-connected Hausdorff spaces, and the subspaces $K\subseteq X$ and $L\subseteq Y$ are compact.  Then we have a commutative diagram
$$\begin{array}{ccc}
{\mathbb K}_\mathrm{hom} (K) & \stackrel{{\cal F}_K\otimes}{\longrightarrow} & {\mathbb K} C^\star_\mathrm{max} \pi (X)|_K \\
\downarrow & & \downarrow \\
{\mathbb K}_\mathrm{hom} (L) & \stackrel{{\cal F}_L\otimes}{\longrightarrow} & {\mathbb K} C^\star_\mathrm{max} \pi (Y)|_L \\
\end{array}$$
\end{lemma}

\begin{proof}
Consider the induced $C^\ast$-functors
$$\begin{array}{ccccc}
1\otimes s_\star & \colon & C(K)\otimes C^\star_\mathrm{max} \pi (X)|_K & \rightarrow & C(K)\otimes C^\star_\mathrm{max} \pi (Y)|_L \\
s^\star\otimes 1 & \colon & C(L)\otimes C^\star_\mathrm{max} \pi (Y)|_L & \rightarrow & C(K)\otimes C^\star_\mathrm{max} \pi (Y)|_L \\
\end{array}$$

Then by proposition \ref{slant2} we have commutative diagrams
$$\begin{array}{ccc}
{\mathbb K}_\mathrm{hom} (K) & \stackrel{{\cal F}_K\otimes}{\longrightarrow} & {\mathbb K} C^\star_\mathrm{max} \pi (X)|_K \\
\| & & \downarrow \\
{\mathbb K}_\mathrm{hom} (K) & \stackrel{(1\otimes s_\star )[{\cal F}_K]\otimes}{\longrightarrow} & {\mathbb K} C^\star_\mathrm{max} \pi (Y)|_L \\
\end{array}$$
and
$$\begin{array}{ccc}
{\mathbb K}_\mathrm{hom} (K) & \stackrel{(s^\star \otimes 1)[{\cal F}_L]\otimes}{\longrightarrow} & {\mathbb K} C^\star_\mathrm{max} \pi (Y)|_L \\
\downarrow & & \downarrow \\
{\mathbb K}_\mathrm{hom} (L) & \stackrel{{\cal F}_L\otimes}{\longrightarrow} & {\mathbb K} C^\star_\mathrm{max} \pi (Y)|_L \\
\end{array}$$

Consider the Hilbert $C(K)\otimes C^\star_\mathrm{max}(\pi (Y)|_L )$-module $\cal E$, where the space ${\cal E}(x)$ consists of all sets of maps
$$\{ \mu_y \colon (\tilde{K},y)\rightarrow \Hom (s(x),s(y))_{C^\star_\mathrm{max} \pi (Y)|_L} \ |\ y\in K \}$$
such that
$$\mu_{y'}(g\tilde{y}) = s_\star (g)\mu_y (\tilde{y})$$
for all paths $\tilde{y}\in (\tilde{K},y)$ and groupoid elements $g\in \Hom(y,y')_{\pi (X)|_K}$.

There are obvious inclusions $i\colon (1\otimes s_\star )[{\cal F}_K]\hookrightarrow {\cal E}$ and $j\colon (s^\star\otimes 1)[{\cal F}_L]\hookrightarrow {\cal E}$.  By proposition \ref{ucp}, for each point $x\in K$ we can find open sets $U\subseteq K$ and $V\subseteq L$ such that $x\in U$, $s(x) \in V$, and
$$p^{-1}[U]\approx U\times C^\star_\mathrm{max} \pi_1 (X,x) \qquad p^{-1}[V]\approx V\times C^\star_\mathrm{max} \pi_1 (X,x)$$

Define modules
$$\begin{array}{rcl}
{\cal E}|_U (x) & = & \{ \phi |_{p^{-1}[U]} \ |\ \phi \in {\cal E}(x) \} \\
{\cal F}_K|_U (x) & = & \{ \phi |_{p^{-1}[U]} \ |\ \phi \in {\cal F}_K(x) \} \\
{\cal F}_L|_V (x) & = & \{ \phi |_{p^{-1}[V]} \ |\ \phi \in {\cal F}_L(x) \} \\
\end{array}$$

Then for all compact sets $A\subseteq U$ and $B\subseteq V$ there are isomorphisms
$$\begin{array}{rcl}
{\cal E}|_A & \cong & C(A)\otimes C^\star_\mathrm{max} \pi (Y)|_L \\
{\cal F}_K|_A & \cong & C(A)\otimes C^\star_\mathrm{max} \pi (X)|_K \\
{\cal F}_L|_B & \cong & C(B)\otimes C^\star_\mathrm{max} \pi (Y)|_L \\
\end{array}$$

Therefore the inclusions $i\colon (1\otimes s_\star )[{\cal F}_K]\hookrightarrow {\cal E}$ and $j\colon (s^\star\otimes 1)[{\cal F}_L]\hookrightarrow {\cal E}$ are isomorphisms.  Hence the Hilbert modules $(1\otimes s_\star )[{\cal F}_K]$ and $(s^\star \otimes 1)[{\cal F}_L]$ are actually equal so we have a commutative diagram
$$\begin{array}{ccc}
{\mathbb K}_\mathrm{hom} (K) & \stackrel{{\cal F}_K\otimes}{\longrightarrow} & {\mathbb K} C^\star_\mathrm{max} \pi (X)|_K \\
\downarrow & & \downarrow \\
{\mathbb K}_\mathrm{hom} (L) & \stackrel{{\cal F}_L\otimes}{\longrightarrow} & {\mathbb K} C^\star_\mathrm{max} \pi (Y)|_L \\
\end{array}$$
as claimed.
\end{proof}

Forming direct limits we obtain a morphism
$$\beta \colon {\mathbb K}_\mathrm{hom}(X)\rightarrow {\mathbb K}C^\star_\mathrm{max}\pi (X)$$

\begin{theorem} \label{ta-da}
The map $\beta$ is natural, and is a weak equivalence when the space $X$ is a single point.  The induced map between stable homotopy groups
$$\beta_\star \colon K_\star (X)\rightarrow K_\star C^\star_\mathrm{max} \pi_1 (X)$$
is the analytic assembly map.
\noproof
\end{theorem}

\begin{proof}
Naturality follows immediately from lemma \ref{biglem2}.  The fact that the induced map between stable homotopy groups is the usual analytic assembly map is a consequence of lemma \ref{biglem1}.  The fact that the map $\beta$ is a weak equivalence when the space $X$ is a single point now follows from the fact that the analytic assembly map $\beta_\star$ is an isomorphism when the space $X$ is a single point. 
\end{proof}

The above theorem can be used to give a characterisation of the analytic assembly map.  The main tool is the following result of Weiss and Williams from \cite{WW}.

\begin{theorem}
Let $F$ be a homotopy-invariant functor from spaces to spectra.\footnote{In this theorem all spaces have the homotopy type of $CW$-complexes}  Then there exists a homotopy-invariant functor $F^\%$ from spaces to spectra and a natural transformation $\alpha \colon F^\% \rightarrow F$ such that

\begin{itemize}

\item[$\bullet$] The functors $X\mapsto \pi_\star F^\% (X)$ form a generalised homology theory.

\item[$\bullet$] The map $\alpha \colon F^\% (\textrm{point})\rightarrow F(\textrm{point})$ is a stable equivalence.

\end{itemize}

Further, the functor $F^\%$ and map $\alpha$ are unique up to stable equivalence.
\noproof
\end{theorem}

The transformation $\alpha$ is called the {\em assembly map} associated to the functor $F$.  We can use the above theorem on the functor $X\mapsto {\mathbb K}C^\star_\mathrm{max} \pi (X)$ together with theorem \ref{ta-da} to obtain the main result of this article.

\begin{theorem}
Let $F^\%$ be a homotopy-invariant functor from spaces to spectra such that the functors $X\mapsto \pi_\star F^\% (X)$ form a generalised homology theory.  Let $\alpha \colon F^\% (X)\rightarrow {\mathbb K}C^\star_\mathrm{max}\pi (X)$ be a natural morphism of spectra which is a stable equivalence when the space $X$ is a single point.  Then the map $\alpha$ is equivalent to the analytic assembly map.
\noproof
\end{theorem}

Roughly speaking this result says that anything which resembles the analytic assembly map at the level of spectra actually is the analytic assembly map.

\bibliographystyle{plain}

\bibliography{data}

\end{document}